\documentclass[preprint]{elsarticle}

\usepackage{lineno,hyperref}
\modulolinenumbers[5]

\usepackage{amsmath}
\usepackage{amsfonts}
\usepackage{empheq}
\usepackage{stmaryrd}
\usepackage{physics}
\usepackage{mathtools}
\usepackage{subcaption}
\usepackage[export]{adjustbox}
\usepackage{yfonts}
\usepackage{mathrsfs}
\usepackage{url}

\journal{\url{www.arXiv.org}}

\begin{document}

\begin{frontmatter}

\title{Discontinuous Galerkin methods for a dispersive wave hydro-morphodynamic model with bed-load transport}

\author[Oden]{Kazbek Kazhyken\corref{corrauth}}
\cortext[corrauth]{Corresponding author}
\ead{kazbek@oden.utexas.edu}

\author[CAMGSD]{Juha Videman}
\author[Oden]{Clint Dawson}

\address[Oden]{Oden Institute for Computational Engineering and Sciences, The University of Texas at Austin, Austin, TX 78712, USA}
\address[CAMGSD]{CAMGSD/Departamento de Matem\'atica, Universidade de Lisboa, Universidade de Lisboa, 1049-001 Lisbon, Portugal}

\begin{abstract}
A dispersive wave hydro-morphodynamic model coupling the Green-Naghdi equations (the hydrodynamic part) with the sediment continuity Exner equation (the morphodynamic part) is presented. Numerical solution algorithms based on discontinuous Galerkin finite element discretizations of the model are proposed. The algorithms include both coupled and decoupled approaches for solving the hydrodynamic and morphodynamic parts simultaneously and separately from each other, respectively. The Strang operator splitting technique is employed to treat the dispersive terms separately, and it provides the ability to ignore the dispersive terms in specified regions, such as surf zones. Algorithms that can handle wetting-drying and detect wave breaking are presented. The numerical solution algorithms are validated with numerical experiments to demonstrate the ability of the algorithms to accurately resolve hydrodynamics of solitary and regular waves, and morphodynamic changes induced by such waves. The results indicate that the model has the potential to be used in studies of coastal morphodynamics driven by dispersive water waves, given that the hydrodynamic part resolves the water motion and dispersive wave effects with sufficient accuracy up to swash zones, and the morphodynamic model can capture the major features of bed erosion and deposition.     
\end{abstract}

\begin{keyword}
Green-Naghdi equations, Exner equation, shallow water flows, dispersive waves, sediment transport, discontinuous Galerkin methods
\end{keyword}

\end{frontmatter}


\section{Introduction}
Coastal areas have a dynamic morphological nature driven by non-trivial interactions between sediment transport, bed morphodynamics, and water waves forced by astronomical tides, winds, and long-wave currents. Changes in coastal morphology caused by natural and anthropogenic forces have the potential to negatively affect coastal infrastructure and environment. For example, structural integrity of piers, levees and other coastal infrastructure can be compromised by excessive erosion of bed due to scouring. Moreover, sediment transport and bed morphodynamic processes play an important role in harbor planning and construction since excessive sediment deposition in a harbor may significantly increase its operating costs due to necessary dredging. Environmental concerns include shoreline and beach erosion that may damage natural habitats of endangered protected species, and effects of sediment transport on contaminants, i.e. sediment deposits may serve as dangerous contaminant sinks or sources depending on surrounding physico-chemical conditions. This evidence suggests that  mathematical modeling of hydro-morphodynamic processes in coastal areas required to forecast sediment transport and bed evolution has clear engineering relevance.

Sediment transport and bed morphodynamic processes are driven by water flow properties, such as the flow velocity and turbulence, which, in turn, are affected by changes in bed surface levels. Therefore, any mathematical modeling of hydro-morphodynamic processes involves a coupling between a hydrodynamic model, that describes water waves and motion, and a sediment transport and bed evolution model, that resolves changes in bed topology driven by sediment erosion, transport, and deposition rates. A widely used variation of mathematical models for hydro-morphodynamic processes is formed by coupling the nonlinear shallow water equations with the sediment continuity Exner equation. Numerical approaches for solving the resulting set of equations can be fully coupled or decoupled, and use structured or unstructured computational grids. The use of unstructured grids can be highly advantageous as they provide the ability of local grid refinement near important bathymetric features and structures. This ability can improve accuracy while maintaining lower computational costs as compared to models that use structured grid methods. Moreover, unstructured grids have better capacity to capture irregular geometries of coastal areas, which is a great advantage over structured grids when hydro-morphodynamic processes are modeled around coastal areas. 

The majority of the numerical schemes developed for the coupled system of the Exner and nonlinear shallow water equations use finite volume methods. Canestrelli \emph{et al.} \cite{canestrelli_etal_2009, canestrelli_etal_2010} have developed a finite volume PRICE-C scheme over unstructured grids for the fully coupled system for the 1D and 2D cases. Castro D\'iaz \emph{et al.} have introduced, for the 1D case in \cite{diaz_etal_2008} and for the 2D case in \cite{diaz_etal_2010}, a finite volume scheme that utilizes the theory of Dal Maso \emph{et al.} \cite{dal_maso_etal_1995} to handle the nonconservative product present in the source term of the fully coupled system. Kozyrakis \emph{et al.} have developed a finite volume scheme over unstructured grids to study coastal hydro-morphodynamics with the fully coupled system in \cite{kozyrakis_etal_2016}. A number of works have been published on the application of wetting-drying processes for the system, e.g. Liang \cite{liang_2011}, Barzgaran \emph{et al.} \cite{barzgaran_etal_2019}, and Rehman and Cho \cite{rehman_cho_2019}. Serrano-Pacheco \emph{et al.} \cite{serrano-pacheco_etal_2012} have developed an upwinding numerical flux finite volume scheme over unstructured grids for both the coupled and decoupled system. High-performance computing applications with graphical processing units (GPUs) have been studied by Garc\'ia‑Navarro \emph{et al.} \cite{garcia-navarro_etal_2019}  with achieved speedup of $O(10^2)$ compared to legacy systems. 

Examples of decoupled approaches over unstructured grids with discontinuous Galerkin methods include Kubatko \emph{et al.} \cite{kubatko_etal_2006} and Izem \emph{et al.} \cite{izem_etal_2012}. A decoupled approach suggests that the nonlinear shallow water equations and the Exner equation are solved separately from each other. In cases where the morphodynamic model has time scales much longer than the hydrodynamic model, updates in bed elevation may be done every $O(10^2)$ time steps of the hydrodynamic model \cite{kubatko_etal_2006, izem_etal_2012}. Although this may provide the opportunity to reduce the amount of computational resources required to run the decoupled model, the method may not be suitable for rapidly evolving beds. In this case, a fully coupled model that solves the hydrodynamic and morphodynamic models simultaneously is more fitting. The resulting coupled model forms a system of hyperbolic nonconservative partial differential equations due to the presence of a nonconservative product in the source term. This fact adds a degree of complexity to the coupled model's numerical solution algorithm. Among examples of discontinuous Galerkin formulations for the coupled nonconservative system are Tassi \emph{et al.} \cite{tassi_etal_2008}, Rhebergen \emph{et al.} \cite{rhebergen_2008}, and Mirabito \emph{et al.} \cite{mirabito_etal_2011}. A major detail of these methods is the special treatment of the nonconservative product term developed using the theory of Dal Maso \emph{et al.} \cite{dal_maso_etal_1995}.

The choice of the nonlinear shallow water equations is popular for a number of reasons: multiple numerical solution algorithms have been developed for these equations (e.g. discontinuous Galerkin implementations in Aizinger and Dawson \cite{aizinger_and_dawson_2002}, Kubatko \emph{et al.} \cite{kubatko_etal_2006_nswe}), a track record of successful application in real world scenarios (e.g. storm surge modeling in Dawson \emph{et al.} \cite{dawson_etal_2011}), the ability of these equations to handle wetting-drying phenomenon that is important for coastal applications (e.g. Bunya \emph{et al.} \cite{bunya_etal_2009}), efficient parallelization strategies (e.g. hybrid MPI+OpenMP, and HPX parallelization in Bremer \emph{et al.} \cite{bremer_etal_2019}), and the ability to capture wave breaking in surf zones. Although the nonlinear shallow water equations provide this multitude of advantages, their lack of ability to capture dispersive wave effects can be a major disadvantage when water wave dynamics must be modeled in areas where wave dispersion is prevalent. An alternative depth-averaged hydrodynamic model that can capture these effects is formed by the Green-Naghdi equations developed in \cite{green_naghdi_1976}. 

The capacity to capture dispersive wave effects comes, however, with a greater analytical and numerical complexity. Among numerical solution algorithms proposed for the Green-Naghdi equations, a few have been based on the Strang operator splitting technique (e.g. Bonneton \emph{et al.} \cite{bonneton_etal_2011}, Samii and Dawson \cite{samii_and_dawson_2018}). In this approach the Green-Naghdi equations are split into two parts: (1) the nonlinear shallow water equations, and (2) the dispersive correction part of the equations. A numerical solution operator for the Green-Naghdi equations is then defined as a successive application of numerical solution operators for these two parts. Although numerical solution algorithms for the two parts  do not have to employ the same discretization method (e.g. Lannes and Marche \cite{lannes_and_marche_2015} use a finite volume method for the first part and a finite difference method for the second part), Duran and Marche \cite{duran_and_marche_2017} use a discontinuous Galerkin method for both parts, and Samii and Dawson \cite{samii_and_dawson_2018} use a hybridized discontinuous Galerkin method to discretize both parts. The operator splitting approach provides a possibility to switch between the nonlinear shallow water equations and the Green-Naghdi equations when modeling water flow dynamics. The switching to the nonlinear shallow water equations can be simply done by not applying the dispersive correction part in areas where the Green-Naghdi equations provide a less accurate model, e.g. in surf zones where wave breaking occurs \cite{bonneton_etal_2011}.

This work aims to introduce a dispersive wave hydro-morphodynamic model by coupling the Green-Naghdi equations with the sediment continuity Exner equation, and to develop numerical solution algorithms for the model. Major motivation for the derivation of this model is its application in a future work to forecast morphodynamic evolution of coastal areas due to dispersive water waves. A significant portion of this work comprises the development of a massively parallel solver that uses the presented numerical solution algorithms. The solver extends a C++ software package \footnote{The software is under development on the date of the publication, and can be accessed at \url{www.github.com/UT-CHG/dgswemv2}. Should there be any questions, comments, or suggestions, please contact the developers through the repository issues page.} developed by Bremer and Kazhyken, and has the capacity to execute numerical simulations of water waves using discontinuous Galerkin discretizations of the nonlinear shallow water and Green-Naghdi equations. 

This paper is organized as follows. In Section 2, the governing equations are presented for the developed mathematical model. The numerical solution algorithms, using discontinuous Galerkin methods over unstructured grids, are introduced in Section 3 both for the decoupled and coupled models. Section 4 presents numerical experiments that are used to demonstrate the ability of the dispersive wave hydro-morphodynamic model to accurately simulate hydrodynamics of solitary and regular waves, and morphodynamic changes induced by such waves. Moreover, the hydrodynamic part of the model is used to simulate water waves in the vicinity of the Faro-Olh\~ao inlet of the Ria Formosa lagoon in Portugal to demonstrate its ability to be used over irregularly shaped domains. Final conclusions are presented in Section 5.

\section{Governing equations}
\begin{figure}
\center
\includegraphics[width=3in]{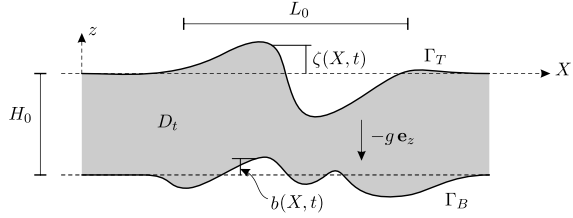}
\caption{A model representation of a body of water as a domain $D_t \subset \mathbb R^{d+1}$.}
\label{Fig:Domain}
\end{figure}

For purposes of this work, a body of water is represented by a domain $D_t \subset \mathbb R^{d+1}$ filled with water as an incompressible, homogeneous, inviscid fluid. In this description, $d$ stands for the horizontal spatial dimension that can take values 1 or 2, $t$ represents the time variable, $\Gamma_T$ and $\Gamma_B$ are the top and bottom boundaries of the domain, respectively, $L_0$ is the characteristic length, and $H_0$ is the reference depth (cf. Fig.\ref{Fig:Domain}). It is assumed that $\Gamma_T$ and $\Gamma_B$ can be represented as graphs, and fluid particles do not cross the boundaries. Both boundaries vary with time: $\Gamma_B$  due to sediment transport and bed morphodynamic processes, $\Gamma_T$  as the evolving free surface of the body of water. The bathymetry, $b(X,t)$, and the free surface elevation, $\zeta(X,t)$, of the body of water are used in the parameterization of $\Gamma_B$ and $\Gamma_T$: 
\begin{linenomath}
\begin{subequations}
\begin{align}
\Gamma_B &= \{(X,-H_0+b(X,t)):X\in \mathbb R^d\}, \\
\Gamma_T &= \{(X,\zeta(X,t)):X\in \mathbb R^d\},
\end{align}
\end{subequations}
\end{linenomath}
and the domain $D_t$ is defined as a set of points $(X,z) \in \mathbb R^d \times \mathbb R$ where $-H_0+b(X,t) < z < \zeta(X,t)$.

Motion of water over an erodible bed and subsequent sediment transport and bed surface evolution are highly interactive processes. Water flow parameters, such as the flow velocity and turbulence, determine the rates of sediment erosion, transport, and deposition that drive changes in bed relief; and these changes, in turn, affect the flow parameters. Therefore, any mathematical modeling of these interrelated hydro-morphodynamic processes involves some sort of coupling between a hydrodynamic model, which governs the changes in flow parameters, and a sediment transport and bed morphodynamic model, which determines the sediment erosion, transport, and deposition rates, and the subsequent changes in bed levels.

\subsection{Hydrodynamic model}
Defining the shallowness parameter $\mu=H_0^2/L_0^2$, the shallow water flow regime is in action when $\mu\ll 1$. Under the assumption of the shallow water flow regime, the Green-Naghdi equations, a depth-averaged hydrodynamic model, provide a sufficiently accurate approximation to water flow dynamics within the domain $D_t$ while maintaining the ability to capture wave dispersion effects \cite{bonneton_etal_2011}. A single parameter variation of the Green-Naghdi equations introduced by Bonneton \emph{et al.} in \cite{bonneton_etal_2011} are defined over a horizontal domain $\Omega \subset \mathbb R^d$ as 
\begin{linenomath}
\begin{equation}\label{Eq:GN}
\partial_t \boldsymbol q + \nabla \cdot \boldsymbol F(\boldsymbol q) + \boldsymbol{D}(\boldsymbol q)  = \boldsymbol S(\boldsymbol q),
\end{equation}
\end{linenomath}
where
\begin{linenomath}
\begin{equation}
\boldsymbol q = \begin{Bmatrix} h \\ h \mathbf u \end{Bmatrix}, \quad
\boldsymbol F(\boldsymbol q) = \begin{Bmatrix} h \mathbf u \\ h\mathbf u \otimes \mathbf u + \frac 1 2 g h^2 \mathbf I \end{Bmatrix}, \quad
\boldsymbol S(\boldsymbol q) = \begin{Bmatrix} 0 \\ -gh \nabla b + \mathbf f \end{Bmatrix},
\end{equation}
\end{linenomath}
$\mathbf u$ is the water velocity represented by a $d$ dimensional vector, $h$ is the water depth represented by the mapping $h(X,t) = \zeta(X,t) + H_0 - b(X,t)$ and assumed to be bounded from below by a positive value, $\mathbf f$ comprises additional source terms for the  momentum continuity equation, e.g. the Coriolis, bottom friction, and surface wind stress forces, $g$ is the acceleration due to gravity, $\mathbf I \in \mathbb R^{d\times d}$ is the identity matrix, and where the wave dispersion effects are introduced into the model through the dispersive term 
\begin{linenomath}
\begin{equation}
\boldsymbol D(\boldsymbol q) = \begin{Bmatrix} 0 \\ \mathbf w_1 - \alpha^{-1} g h \nabla \zeta \end{Bmatrix}.
\end{equation}
\end{linenomath}
In this description, $\mathbf w_1$ is defined through
\begin{linenomath}
\begin{equation}\label{Eq:w1}
(\mathbf I + \alpha h \mathcal T h^{-1}) \mathbf w_1 = \alpha^{-1} g h \nabla \zeta + h \mathcal Q_1(\mathbf u),
\end{equation}
\end{linenomath}
where operators $\mathcal T$ and $\mathcal Q_1$ are
\begin{linenomath}
\begin{subequations}
\begin{align}
\mathcal T(\mathbf w) =&\,\mathcal R_1(\nabla\cdot\mathbf w) + \mathcal R_2(\nabla b \cdot \mathbf w), \\
\mathcal Q_1(\mathbf w) =&-2\mathcal R_1\left(\partial_{x} \mathbf w \cdot \partial_{y} \mathbf w^\perp+(\nabla \cdot \mathbf w)^2\right)+\mathcal R_2\left(\mathbf w\cdot (\mathbf w \cdot \nabla)\nabla b\right),
\end{align}
\end{subequations}
\end{linenomath}
with operators $\mathcal R_1$ and $\mathcal R_2$ defined as
\begin{linenomath}
\begin{subequations}
\begin{align}
\mathcal R_1(w) &= -\frac 1 {3h} \nabla(h^3 w) -  \frac {h}{2} w \nabla b,\\
\mathcal R_2(w) &= \frac 1 {2h} \nabla (h^2 w) +  w \nabla b.
\end{align}
\end{subequations}
\end{linenomath}
In this description, $\alpha\in\mathbb R$ is a parameter that is used to optimize dispersive properties of the presented hydrodynamic model. By adjusting $\alpha$, the difference between the phase and group velocities coming from the Stokes linear theory and the Green-Naghdi equations can be minimized. Minimizing the averaged variation over some range of wave number values is a common strategy that puts practical values for $\alpha$ in the range between $1.0$ and $1.2$ \cite{bonneton_etal_2011}.

\subsection{Sediment transport and bed morphodynamic model}
Among modes of sediment transport are bed-load, suspended-load, and wash-load transport. In the presented work, the developed model is limited to bed-load transport, where sediment particles slide, roll, and saltate due to shearing forces from the surrounding fluid while staying sufficiently close to bed. The sediment continuity Exner equation provides a mathematical model that describes morphological evolution of bed due to sediment transport phenomena \cite{exner_1925}. In a morphodynamic model limited to bed-load transport, the equation states that change of $b(X,t)$ in time is equal to the divergence of the sediment flux $\mathbf Q_b$:
\begin{linenomath}
\begin{equation}\label{Eq:Exner}
\partial_{t} b + \nabla\cdot \mathbf Q_b = 0,
\end{equation}
\end{linenomath}
where $\mathbf Q_b$ is an empirically defined function \cite{diaz_etal_2008}. Intuitively, sediment transport occurs in the flow direction; therefore, 
\begin{linenomath}
\begin{equation}\label{Eq:SedFlux}
\mathbf Q_b = \lvert \mathbf Q_b \rvert \mathbf {\bar u},
\end{equation}
\end{linenomath}
where $\mathbf {\bar u}$ is the unit flow velocity vector, and $\lvert \mathbf Q_b \rvert$ is the magnitude of the sediment flux represented by an empirical formula. A number of empirical models have been proposed for $\lvert \mathbf Q_b \rvert$; most of them may be represented as (see \cite{diaz_etal_2008, cordier_etal_2011} and all references therein)
\begin{linenomath}
\begin{equation}\label{Eq:Qb}
\lvert \mathbf Q_b \rvert = A(h, \mathbf u)\lvert \mathbf u \rvert^{m},
\end{equation}
\end{linenomath}
where $1 \leq m \leq 3$ and $A(h, \mathbf u)$ is an empirical equation, e.g. the Grass model takes $A$ as a constant calibrated for the application under investigation and sets $m=3$ \cite{grass_1981}. There are a number of other empirical expressions for $\lvert \mathbf Q_b \rvert$, e.g. Meyer-Peter and Mueller \cite{meyer_and_muller_1948}, Fernandez Luque and Van Beek \cite{luque_beek_1976}, Nielsen \cite{nielsen_1992}, Ribberink \cite{ribberink_1998}. The choice of the empirical representation of $\lvert \mathbf Q_b \rvert$ is judicious and influenced by the application.

\section{Numerical methods}
Discontinuous Galerkin finite element methods are used for discretizing the governing equations. This choice facilitates the use of unstructured meshes that are well suited for irregular geometries of coastal areas. 

\subsection{Notation and functional setting}
The problem domain $\Omega$ is partitioned into a finite element mesh $\mathcal T_h = \{K\}$ that provides an approximation to the domain: 
\begin{linenomath}
\begin{equation}
\Omega\approx\Omega_h=\sum_{K\in \mathcal{T}_h}K,
\end{equation}
\end{linenomath}
where the subscript $h$ stands for the mesh parameter represented by the diameter of the smallest element in the mesh. The set of all faces of elements of the mesh, $\partial\mathcal T_h$, and the set of all edges of the mesh skeleton, $\mathcal{E}_h$, are defined as
\begin{linenomath}
\begin{subequations}
\begin{align}
\partial\mathcal T_h &= \lbrace\partial K : K\in \mathcal{T}_h\rbrace,\\
\mathcal{E}_h &= \lbrace e\in\bigcup_{K\in\mathcal{T}_h}\partial K\rbrace.
\end{align}
\end{subequations}
\end{linenomath}
Note that in $\mathcal E_h$ the common element faces appear only once but in $\partial \mathcal T_h$ they are counted twice. 

To develop variational formulations of the governing equations, inner products are defined for finite dimensional vectors $\boldsymbol u$ and $\boldsymbol v$ through:
\begin{linenomath}
\begin{subequations}
\begin{align}
(\boldsymbol u,\boldsymbol v)_\Omega &= \int_\Omega \boldsymbol u \cdot \boldsymbol v \, \dd X, \\
\langle \boldsymbol u, \boldsymbol v \rangle_{\partial \Omega} &= \int_{\partial \Omega}\boldsymbol u \cdot \boldsymbol v\, \dd X,
\end{align}
\end{subequations}
\end{linenomath}
for $\Omega \subset \mathbb R^d$ and $\partial \Omega \subset \mathbb R^{d-1}$.

An approximating space of trial and test functions is chosen as the set of square integrable functions over $\Omega_h$ such that their restriction to an element $K$ of the mesh belongs to $\mathcal Q^p(K)$, a space of polynomials of degree at most $p \ge 0$ with support in $K$:
\begin{linenomath}
\begin{equation}
\mathbf V_h^{p,m} \coloneqq \{\boldsymbol v \in (L^2(\Omega_h))^{m}: \boldsymbol v|_K \in \mathcal (\mathcal Q^p(K))^{m} \quad  \forall K \in \mathcal T_h \},
\end{equation}
\end{linenomath}
and, similarly, an approximation space over the mesh skeleton is chosen as
\begin{linenomath}
\begin{equation}
\mathbf M_h^{p,m} \coloneqq \{\boldsymbol \mu \in (L^2(\mathcal E_h))^{m}: \boldsymbol \mu|_e \in \mathcal (\mathcal Q^p(e))^{m} \quad \forall e \in \mathcal E_h\}.
\end{equation}
\end{linenomath}

\subsection{Decoupled model}
In the decoupled model method the Green-Naghdi and Exner equations are solved separately. After the flow parameters are evolved in time according to the hydrodynamic model for a number of time steps, the bed surface elevation is updated with the use of the morphodynamic model and fed back into the hydrodynamic model to continue the evolution of the flow parameters until the next bed surface elevation update. In cases where the time scales in the hydrodynamic model are much shorter than the time scales in the morphodynamic model, the bed surface elevation does not need to be updated every time step of the hydrodynamic model. In some cases the bed update may happen every $O(10^2)$ time steps of the hydrodynamic model \cite{kubatko_etal_2006, izem_etal_2012}. The ability to save computational resources is the main advantage of the decoupled model method. However, this method may be unsuitable if the time scales in the hydrodynamic and morphodynamic models are comparable, e.g. in the case of a dam break \cite{kubatko_etal_2006, izem_etal_2012}.

The Green-Naghdi equations presented in Eq.(\ref{Eq:GN}) can be treated numerically with the use of the well-known Strang operator splitting technique \cite{bonneton_etal_2011, samii_and_dawson_2018}. The equation is split into: (1) the nonlinear shallow water equations by dropping the dispersive term of the equation, and (2) the dispersive correction part where the wave dispersion effects on flow velocities are introduced into the model through the dispersive term. If $\mathcal S_1$ is a numerical solution operator for the nonlinear shallow water equations, i.e. $\mathcal S_1(\Delta t)$ propagates numerical solution by a time step $\Delta t$, and, similarly, $\mathcal S_2$ is a numerical solution operator for the dispersive correction part, then the second-order Strang operator splitting technique \cite{strang_1968} states that a numerical solution operator for the Green-Naghdi equations can be approximated as
\begin{linenomath}
\begin{equation}
\mathcal S(\Delta t) = \mathcal S_1(\Delta t/2) \mathcal S_2(\Delta t) \mathcal S_1(\Delta t/2),
\end{equation}
\end{linenomath}
where $\mathcal S$ is a second-order temporal discretization if both $\mathcal S_1$ and $\mathcal S_2$ use a second-order time discretization method.

A numerical solution operator $\mathcal S_1$ for the nonlinear shallow water equations is developed using a discontinuous Galerkin finite element formulation. Therefore, an approximate solution $\boldsymbol q_h \in \mathbf V_h^{p,d+1}$ must satisfy the variational formulation
\begin{linenomath}
\begin{equation}\label{Eq:NSWEVarLoc}
(\partial_t \boldsymbol{q}_h,\boldsymbol{v})_{\mathcal{T}_h}-(\boldsymbol{F}_h,\nabla \boldsymbol{v})_{\mathcal{T}_h}+\langle\boldsymbol{F}_h^*,\boldsymbol{v}\rangle_{\partial {\mathcal{T}_h}}-(\boldsymbol{S}_h,\boldsymbol{v})_{\mathcal{T}_h}=0\quad \forall \boldsymbol{v}\in\textbf{V}_h^{p,d+1},
\end{equation}
\end{linenomath}
where $\boldsymbol{F}_h=\boldsymbol F(\boldsymbol{q}_h)$ and $\boldsymbol{S}_h=\boldsymbol S(\boldsymbol{q}_h)$, ${\boldsymbol F_h^*}$ is a single valued approximation to $\boldsymbol F_h \mathbf n$ over element faces, called the numerical flux, and $\mathbf n$ is the unit outward normal vector to element face. The present work uses the numerical flux from the hybridized discontinuous Galerkin method developed by Samii \emph{et al.} in \cite{samii_etal_2019}. Therefore, the numerical flux is defined through $\widehat {\boldsymbol q}_h\in \mathbf M_h^{p,d+1}$, an approximation to $\boldsymbol q$ over the mesh skeleton called the numerical trace, as in \cite{samii_etal_2019}
\begin{linenomath}
\begin{equation}
\boldsymbol F_h^* = \widehat{\boldsymbol F}_h\mathbf n + \boldsymbol \tau (\boldsymbol q_h - \widehat{\boldsymbol q}_h),
\end{equation}
\end{linenomath}
where $\widehat{\boldsymbol F}_h=\boldsymbol F(\widehat{\boldsymbol q}_h)$, and $\boldsymbol \tau$ is a stabilization parameter motivated by the local Lax-Friedrichs numerical flux:
\begin{linenomath}
\begin{equation}
\boldsymbol \tau =\lambda_{\max}(\widehat{\boldsymbol q}_h).
\end{equation}
\end{linenomath}
In this description of the stabilization parameter, $\lambda_{\max}$ is the maximum eigenvalue of the normal Jacobian matrix $\boldsymbol A = \partial_{\boldsymbol q} (\boldsymbol F \mathbf n)$:
\begin{linenomath}
\begin{equation}
\lambda_{\max}(\boldsymbol q) = \lvert\mathbf{u}\cdot\mathbf{n}\rvert+\sqrt{gh}.
\end{equation}
\end{linenomath}
The numerical trace $\widehat {\boldsymbol q}_h\in \mathbf M_h^{p,d+1}$ must be such that the numerical flux is conserved across all internal edges in the mesh skeleton, and boundary conditions are satisfied at all boundary edges through the boundary operator $\boldsymbol B_h$ defined according to an imposed boundary condition as in \cite{samii_etal_2019}: 
\begin{linenomath}
\begin{equation}\label{Eq:NSWEVarGlob}
\langle \boldsymbol F^*_h , \boldsymbol \mu \rangle_{\partial \mathcal T_h \backslash \partial \Omega_h} + \langle \boldsymbol B_h , \boldsymbol \mu \rangle_{\partial \mathcal T_h \cap \partial \Omega_h}=0\,\,\,\forall \boldsymbol\mu\in \mathbf M_h^{p,d+1}.
\end{equation}
\end{linenomath}
Eq.(\ref{Eq:NSWEVarLoc}) and Eq.(\ref{Eq:NSWEVarGlob}) form a system of equations that is used to solve for an approximate solution $\boldsymbol q_h \in \mathbf V_h^{p,d+1}$. For complete details of the formulation along with definitions for $\boldsymbol B_h$, see Samii \emph{et al.} \cite{samii_etal_2019}.

In order to generate $\mathcal S_2$, a numerical solution operator for the dispersive correction part of the Green-Naghdi equations, Eq.(\ref{Eq:w1}) is written as a system of first order equations using the definition for operator $\mathcal T$ \cite{samii_and_dawson_2018}:
\begin{linenomath}
\begin{empheq}[left=\empheqlbrace]{equation}\label{Eq:w1w2}
\begin{split}
&\nabla \cdot (h^{-1} \mathbf w_1) - h^{-3} w_2 = 0\\
&\mathbf w_1- \tfrac 1 3 \nabla w_2 - \tfrac {1}{2} h^{-1} w_2 \nabla b 
+ \tfrac 1 2 \nabla (h \nabla b \cdot \mathbf w_1) +
 \mathbf w_1 \nabla b \otimes \nabla b = \mathbf{s}(\boldsymbol{q})
\end{split},
\end{empheq}
\end{linenomath}
where $\textbf{s}(\boldsymbol{q}) = \alpha^{-1}gh \nabla \zeta + h \mathcal Q_1(\mathbf u)$. A variational formulation for Eq.(\ref{Eq:w1w2}) forms a global system of equations that would benefit from a dimensional reduction. Therefore, the hybridized discontinuous Galerkin method developed by Samii and Dawson in  \cite{samii_and_dawson_2018} is employed to treat numerically Eq.(\ref{Eq:w1w2}). According to \cite{samii_and_dawson_2018}, an approximate solution $(\mathbf w_{1h}, w_{2h})\in \mathbf V_h^{p,d+1}$ and $\hat{\mathbf w}_{1h} \in {\mathbf M}_h^{p,d}$ are sought such that
\begin{linenomath}
\begin{empheq}[left=\empheqlbrace]{equation}\label{Eq:GNloc}
\begin{split}
&(h^{-3}\, w_{2h},v_2)_{\mathcal T_h} - \langle{\hat h}^{-1}\, \hat{\mathbf w}_{1h} \cdot \mathbf n, v_2\rangle_{\partial \mathcal T_h}+ \left(h^{-1}\, \mathbf w_{1h}, \nabla v_2\right)_{\mathcal T_h} = 0 \\
&\begin{aligned}
(\mathbf{w}_{1h},\mathbf v_{1})_{\mathcal T_h} &- \langle{\tfrac 1 3 \mathbf w^*_{2h}} , \mathbf v_1 \rangle_{\partial \mathcal T_h} + \\
&+\left(\tfrac 1 3 w_{2h},\nabla \cdot \mathbf v_1\right)_{\mathcal T_h}-\left( \tfrac 1 2 h^{-1} \nabla b\, w_{2h}, \mathbf v_1\right)_{\mathcal T_h} + \\
&+\langle\tfrac 1 2 \hat h \nabla b \cdot \hat{\mathbf w}_{1h}, \mathbf v_1 \cdot \mathbf n \rangle_{\partial \mathcal T_h}-\left(\tfrac 1 2 h \nabla b \cdot \mathbf w_{1h},\nabla \cdot \mathbf v_1\right)_{\mathcal T_h} +\\
&+ \left(\nabla b \otimes \nabla b\,\mathbf{w}_{1h} , \mathbf v_1 \right)_{\mathcal T_h}=\left(\mathbf{s}_h, \mathbf v_1\right)_{\mathcal T_h}
\end{aligned}
\end{split},
\end{empheq}
\end{linenomath}
for all $(\mathbf v_1, v_2) \in \mathbf V_h^{p,d+1}$, where $\mathbf{s}_h=\mathbf{s}(\boldsymbol{q}_h)$, and the numerical flux $\mathbf{w}_{2h}^* $ is defined as
\begin{linenomath}
\begin{equation}
\mathbf w_{2h}^* = w_{2h}\mathbf n + \boldsymbol \tau \left(\mathbf{w}_{1h} - \hat{\mathbf{w}}_{1h}\right),
\end{equation}
\end{linenomath}
with a scalar constant $\boldsymbol{\tau}$ used as the stabilization parameter. The numerical flux is weakly conserved and the imposed boundary conditions, defined through the boundary operator $\mathcal B_h$, are weakly satisfied as in \cite{samii_and_dawson_2018}:
\begin{linenomath}
\begin{equation}\label{Eq:GNglob}
\langle \mathbf w_{2h}^* , \boldsymbol \mu \rangle_{\partial \mathcal T_h\backslash \partial \Omega_h} + \langle\mathcal B_h,\boldsymbol \mu\rangle_{\partial \mathcal T_h \cap \partial \Omega_h} = 0 \,\,\,\forall \boldsymbol \mu\in {\mathbf M}_h^{p,d}.
\end{equation}
\end{linenomath}
Eq.(\ref{Eq:GNloc}) is a series of local systems which forms block diagonal matrices that can be used to perform efficient static condensation of Eq.(\ref{Eq:GNglob}). This will form a global system of equations with its dimension equal to the dimension of ${\mathbf M}_h^{p,d}$. The system is solved to obtain $\hat{\mathbf w}_{1h} \in {\mathbf M}_h^{p,d}$ that is subsequently substituted back into Eq.(\ref{Eq:GNloc}) to recover $\mathbf{w}_{1h}\in\textbf{V}_h^{p,d}$. The result is then used in the dispersive correction portion of the Green-Naghdi equations to seek an approximate solution $\boldsymbol{q}_h\in\textbf{V}_h^{p,d+1}$ that satisfies the variational formulation
\begin{linenomath}
\begin{equation}
(\partial_t \boldsymbol{q}_h,\boldsymbol{v})_{\mathcal{T}_h}+\left(\boldsymbol D_h,\boldsymbol{v}\right)_{\mathcal{T}_h}=0\quad \forall \boldsymbol{v}\in\textbf{V}_h^{p,d+1},
\end{equation}
\end{linenomath}
where $\boldsymbol D_h = \boldsymbol D(\boldsymbol{q}_h)$. High order derivatives of $\mathbf{u}_h$, present in $\mathcal Q_1(\mathbf{u}_h)$, are computed weakly using a discontinuous Galerkin method with centered numerical fluxes. See \cite{samii_and_dawson_2018} for complete details of the presented formulation along with definitions of the boundary operators $\mathcal B_h$.

As a scalar conservation law, the Exner equation can be efficiently discretized using a discontinuous Galerkin method. To this end, an approximate solution $b_h\in{\mathbf V}_h^{p,1}$ is sought such that 
\begin{linenomath}
\begin{equation}\label{Eq:ExnerVar}
\left(\partial_{t} b_h,v\right)_{\mathcal T_h} -\left(\mathbf Q_b, \nabla v\right)_{\mathcal T_h}+ \langle\mathbf Q_b^*,v\rangle_{\partial \mathcal T_h} = 0\,\,\,\forall v\in {\mathbf V}_h^{p,1},
\end{equation}
\end{linenomath}
where  a simple upwinding scheme is employed for the numerical flux $\mathbf Q_b^*$ since the sediment flux is not an explicit function of $b$ and the normal Jacobian matrix cannot be formed. Since the bed-load transport models in Eq.(\ref{Eq:Qb}) move sediment particles in the flow direction, the numerical flux $\mathbf Q_b^*$ is defined as \cite{mirabito_etal_2011}:
\begin{linenomath}
\begin{equation}
\mathbf Q_b^*=
\begin{cases}
\mathbf Q_b^+&\text{if}\,\,\,\,\mathbf {\hat u} \cdot \mathbf n \geq 0 \\
\mathbf Q_b^-&\text{if}\,\,\,\,\mathbf {\hat u} \cdot \mathbf n < 0 
\end{cases},
\end{equation}
\end{linenomath}
where $\mathbf {\hat u}$ is the Roe-averaged velocity defined as
\begin{linenomath}
\begin{equation}
\mathbf {\hat u} = \frac{\mathbf{u}^+\sqrt{h^+} + \mathbf{u}^-\sqrt{h^-}}{\sqrt{h^+}+\sqrt{h^-}}.
\end{equation}
\end{linenomath}
In this description, superscript $+$ denotes a variable value at $\partial K$ when approaching from the interior of an element $K$, and $-$ when approaching from the exterior.

\subsection{Coupled model}
In the coupled model method, the Green-Naghdi and Exner equations, Eq.(\ref{Eq:GN}) and Eq.(\ref{Eq:Exner}), are fully coupled and solved simultaneously. The Strang operator splitting technique is used also for the coupled model and the numerical solution operator for the dispersive correction part, $\mathcal S_2$, is as in the decoupled model. However, the operator $\mathcal S_1$ has to be modified since it now needs to provide a numerical solution to the coupled system of the nonlinear shallow water and Exner equations and not only to the nonlinear shallow water equations as in the decoupled model. 

The discontinuous Galerkin method developed for hyperbolic nonconservative partial differential equations  by Rhebergen \emph{et al.} \cite{rhebergen_2008} is used, in the form presented by Mirabito \emph{et al.} \cite{mirabito_etal_2011}, for the model that couples the nonlinear shallow water and Exner equations. In this method, the numerical scheme for the Exner equation is as in the  decoupled model but the numerical scheme for the nonlinear shallow water equations requires corrections due to the nonconservative term $-gh \nabla b$ present in the source term $\boldsymbol{S}(\boldsymbol q)$. Defining $\boldsymbol{p} = \{\boldsymbol q\quad b\}^T$, introducing a third order tensor $\boldsymbol{G}(\boldsymbol{p})$ such that $\boldsymbol G(\boldsymbol p) \nabla \boldsymbol{p} = \{0\quad-gh \nabla b\}^T$, and setting $\boldsymbol{s} = \{0\quad\mathbf{f}\}^T$, we require that an approximate solution $\boldsymbol{q}_h\in\textbf{V}_h^{p,d+1}$ to the nonlinear shallow water equations satisfies the variational formulation \cite{mirabito_etal_2011}
\begin{linenomath}
\begin{equation}
\begin{aligned}
\left(\partial_t \boldsymbol{q}_h, \boldsymbol{v}\right)_{\mathcal{T}_h}-&\left(\boldsymbol{F}_h,\nabla \boldsymbol{v}\right)_{\mathcal{T}_h}+\left\langle\boldsymbol{F}_h^*,\boldsymbol{v}\right\rangle_{\partial {\mathcal{T}_h}}+\\
-&\langle\smallint_0^1\boldsymbol{G}(\phi(\tau;\boldsymbol{p}_h^L,\boldsymbol{p}_h^R))\frac{d\phi}{d\tau}(\tau;\boldsymbol{p}_h^L,\boldsymbol{p}_h^R)\dd\tau\,\mathbf{n}^L,\boldsymbol{v}^*\rangle_{\mathcal{E}_h\backslash \partial \Omega_h }+\\
-&\left(\boldsymbol G_h \nabla \boldsymbol{p}_h,\boldsymbol{v}\right)_{\mathcal{T}_h}-\left(\boldsymbol{s},\boldsymbol{v}\right)_{\mathcal{T}_h}=0\quad \forall \boldsymbol{v}\in\textbf{V}_h^{p,d+1},
\end{aligned}
\end{equation}
\end{linenomath}
where $\boldsymbol G_h = \boldsymbol G(\boldsymbol{p}_h)$, $\phi(\tau;\boldsymbol{p}_h^L,\boldsymbol{p}_h^R)$ is a Lipschitz continuous path from $\boldsymbol{p}^{L}_h$ to $\boldsymbol{p}^{R}_h$ such that $\phi(0)=\boldsymbol{p}^{L}_h$ and $\phi(1)=\boldsymbol{p}^{R}_h$, and where $\boldsymbol{v}^*=\frac{1}{2}(\boldsymbol{v}^L+\boldsymbol{v}^R)$ with the superscripts $L$ and $R$ corresponding to elements $K^L$ and $K^R$ such that $e =\partial K^L \cap \partial K^R$. The choice of the form for the path $\phi(\tau;\boldsymbol{p}_h^L,\boldsymbol{p}_h^R)$ has minor effect on numerical solutions \cite{rhebergen_2008}; therefore, a simple linear path $\phi(\tau;\boldsymbol{p}_h^L,\boldsymbol{p}_h^R)=(1-\tau)\boldsymbol{p}_h^L+\tau\boldsymbol{p}_h^R$ has been chosen for this numerical formulation. Subsequently, the integral in the nonconservative term may be evaluated as \cite{mirabito_etal_2011}
\begin{linenomath}
\begin{equation}
\begin{aligned}
\boldsymbol{w}_{nc}=\int_0^1\boldsymbol{G}(\phi(\tau;\boldsymbol{p}_h^L,\boldsymbol{p}_h^R)) & \frac{d\phi}{d\tau}(\tau;\boldsymbol{p}_h^L,\boldsymbol{p}_h^R)\dd\tau\,\mathbf{n}^L=\\
&=\begin{Bmatrix} 0 \\ \frac{1}{2}g(h^L+h^R)(b^L-b^R)\mathbf{n}^L\end{Bmatrix}.
\end{aligned}
\end{equation}
\end{linenomath}
It is worth noting that $\boldsymbol{w}_{nc}$ is single valued over the edges of the mesh skeleton and does not depend on the way the elements $K^L$ and $K^R$ are chosen for the edge $e =\partial K^L \cap \partial K^R$. The numerical flux $\boldsymbol{F}_h^*$ for this numerical scheme is defined as \cite{rhebergen_2008}
\begin{linenomath}
\begin{equation}
\boldsymbol{F}_h^* = \begin{cases}
\boldsymbol{F}_h^+\mathbf{n}-\frac{1}{2}\boldsymbol{w}_{nc}&\text{if}\,\,\,\,S^+>0\\
\boldsymbol{F}_h^{\text{HLL}}-\frac{S^++S^-}{2(S^--S^+)}\boldsymbol{w}_{nc}&\text{if}\,\,\,\,S^+\leq0\leq S^-\\
\boldsymbol{F}_h^-\mathbf{n}+\frac{1}{2}\boldsymbol{w}_{nc}&\text{if}\,\,\,\,S^-<0
\end{cases},
\end{equation}
\end{linenomath}
where the truncated characteristic speeds $S^+$ and $S^-$ are 
\begin{linenomath}
\begin{subequations}
\begin{align}
S^+&=\min(\mathbf{u}^+\cdot\mathbf{n}-\sqrt{gh^+}, \mathbf{u}^-\cdot\mathbf{n}-\sqrt{gh^-}),\\
S^-&=\max(\mathbf{u}^+\cdot\mathbf{n}+\sqrt{gh^+}, \mathbf{u}^-\cdot\mathbf{n}+\sqrt{gh^-}),
\end{align}
\end{subequations}
\end{linenomath}
and the Harten–Lax–van Leer flux $\boldsymbol{F}_h^{\text{HLL}}$ is \cite{harten_etal_1983}
\begin{linenomath}
\begin{equation}
\boldsymbol{F}_h^{\text{HLL}}=\frac{1}{S^--S^+}((S^-\boldsymbol{F}_h^+-S^+\boldsymbol{F}_h^-)\mathbf{n}-S^+S^-(\boldsymbol{q}_h^+-\boldsymbol{q}_h^-)).
\end{equation}
\end{linenomath}
Finally, the modified numerical solution operator $\mathcal S_1$ seeks an approximate solution $(\boldsymbol{q}_h , b_h)\in\textbf{V}_h^{p,d+2}$ such that \cite{mirabito_etal_2011}
\begin{linenomath}
\begin{equation}
\begin{aligned}
\left(\partial_t \begin{Bmatrix} \boldsymbol{q}_h \\ b_h\end{Bmatrix}, \boldsymbol{v}\right)_{\mathcal{T}_h}-&\left(\begin{Bmatrix}\boldsymbol{F}_h \\ \mathbf{Q}_b\end{Bmatrix},\nabla \boldsymbol{v}\right)_{\mathcal{T}_h}+\left\langle\begin{Bmatrix}\boldsymbol{F}_h^* \\ \mathbf{Q}_b^*\end{Bmatrix},\boldsymbol{v}\right\rangle_{\partial {\mathcal{T}_h}}+\\
-&\left\langle\begin{Bmatrix}\boldsymbol{w}_{nc} \\ 0\end{Bmatrix},\boldsymbol{v}^*\right\rangle_{\mathcal{E}_h\backslash \partial \Omega_h }-\left(\begin{Bmatrix} \boldsymbol G_h \nabla \boldsymbol{p}_h \\ 0\end{Bmatrix},\boldsymbol{v}\right)_{\mathcal{T}_h}+\\
-&\left(\begin{Bmatrix}\boldsymbol{s} \\ 0\end{Bmatrix},\boldsymbol{v}\right)_{\mathcal{T}_h}=0\quad \forall \boldsymbol{v}\in\textbf{V}_h^{p,d+2}.
\end{aligned}
\end{equation}
\end{linenomath}

\subsection{Wetting-drying, wave breaking, and slope limiting}
In the developed hydro-morphodynamic model,  the water depth $h$ is assumed to be bounded from below by a positive value in the Green-Naghdi equations. This assumption must be ensured with a wetting-drying algorithm that preserves the positivity criterion for the water depth. Since the water depth is updated in the numerical solution operator $\mathcal S_1$ only, the algorithm shall be executed in conjunction with the operator $\mathcal S_1$. In the presented work, the wetting-drying algorithm developed by Bunya \emph{et al.} in \cite{bunya_etal_2009} for discontinuous solutions to the nonlinear shallow water equations is adopted. Among the main features of the algorithm are: (1) the water depth is never allowed to drop below a specified minimum water depth $h_0$, (2) the elements of the mesh used for numerical simulations are defined as "wet" or "dry" according to a classification algorithm, (3) the water mass is allowed to transfer from "wet" to "dry" elements only; otherwise, interfaces between "wet" and "dry" elements are treated as a reflecting boundary. For the dispersive correction and Exner equations a positivity preserving wetting-drying algorithm is not required, and the wetting-drying fronts are modeled as reflecting boundaries.    

Although the Green-Naghdi equations have the ability to capture dispersive properties of water waves, the equations do not accurately resolve wave breaking phenomena in surf zones \cite{bonneton_etal_2011}. A more suitable depth-averaged hydrodynamic model capable of capturing wave breaking phenomena is formed by the nonlinear shallow water equations: as a hyperbolic system, near wave breaking points numerical solutions to these equations develop shocks that provide a sufficiently accurate model of energy dissipation during wave breaking \cite{bonneton_etal_2011}. The use of the Strang operator splitting technique for the numerical treatment of the presented model provides an opportunity to switch between the Green-Naghdi and nonlinear shallow water equations in areas where one model is deemed to be more accurate than the other. In the developed splitting technique, it is possible to switch to the nonlinear shallow water equations by setting $\mathcal S_2=1$ in regions where the Green-Naghdi equations can no longer provide an adequate approximation, e.g. in wave breaking regions. Therefore, a wave breaking detection criterion should be considered. To this end, the wave breaking criterion adopted by Duran and Marche in \cite{duran_and_marche_2017} from the discontinuity detection criterion of Krivodonova \emph{et al.} \cite{krivodonova_etal_2004} is incorporated into the numerical model. The criterion states that wave breaking occurs over an element $K$ if the parameter \cite{duran_and_marche_2017}
\begin{linenomath}
\begin{equation}
\mathbb I_K = \frac{\sum_{F\in\partial K_{\text{in}}}\vert\int_{F}(h^+-h^-)\dd X\vert}{\mathfrak{h}_K^{\frac{p+1}{2}}\,\vert \partial K_{\text{in}}\vert\,\Vert h\Vert_{L^{\infty}(K) }}
\end{equation}
\end{linenomath}
is greater than a specified threshold that is typically $O(1)$. In this description of the parameter $\mathbb I_K$, $\mathfrak{h}_K$ is the element diameter,  $\partial K_{\text{in}}$ are the inflow faces of the element where $\mathbf{u}\cdot\mathbf{n}<0$, and $\vert \partial K_{\text{in}}\vert$ is the total length of the inflow faces. 

In applications of discontinuous Galerkin methods for the nonlinear shallow water equations, a slope limiter may be required in order to remove oscillations at sharp discontinuities in numerical solutions and preserve numerical stability. In particular, the wave breaking phenomena present themselves as sharp discontinuities in the numerical solutions. Therefore, the Cockburn-Shu limiter \cite{cockburn_shu_2001} is incorporated into the numerical model and applied in conjunction with the operator $\mathcal S_1$. Changes in bed elevation may also form shocks that require a limiting procedure to avoid spurious oscillations in numerical solutions; thus, the Xu et al. limiter \cite{xu_etal_2009} is integrated into the model to perform slope limiting in  the Exner equation. The details of the limiters are not presented here, but readers are encouraged to consult the original sources.    

\section{Numerical experiments and discussion}

The developed numerical model has been implemented in a software framework written in C++ programming language with the use of open source scientific computing libraries, such as Eigen \cite{eigen}, Blaze \cite{blaze}, and PETSc \cite{petsc}. The software has been parallelized for shared and distributed memory systems with the use of a hybrid OpenMP+MPI programming, and HPX \cite{hpx}. Performance comparison between the hybrid programming and HPX has been performed by Bremer \emph{et al.} in \cite{bremer_etal_2019}. The software has the capacity to simulate water waves using the discontinuous Galerkin finite element discretizations of the nonlinear shallow water and Green-Naghdi equations  developed in \cite{kubatko_etal_2006_nswe, samii_and_dawson_2018, samii_etal_2019}, and it has been extended with the developed coupled and decoupled numerical models to allow for the possibility to simulate hydro-morphodynamic processes in coastal regions under the action of highly dispersive water waves. 

The decoupled model employs the numerical solution algorithm for the Green-Naghdi equations developed and extensively validated in \cite{samii_and_dawson_2018}; and the coupled model with a rigid bed reduces to the Green-Naghdi equations where the nonlinear shallow water equations are discretized with a discontinuous Galerkin method with the Harten-Lax-van Leer numerical flux, and the dispersive correction part is discretized with an explicit hybridized discontinuous Galerkin method. The latter numerical solution algorithm for the Green-Naghdi equations is validated in the first three examples where simulations are performed over rigid beds and no morphological changes are present. The full dispersive wave hydro-morphodynamic model is validated in the fourth experiment, and the ability of the hydrodynamic model to simulate water waves over irregularly shaped domains is presented in the fifth experiment. In all of the following numerical experiments the Dubiner polynomials of order $p=1$ from \cite{dubiner_1991} are used for the approximating space $\mathbf V_h$, and the Legendre polynomials of order $p=1$ are used for the approximating space $\mathbf M_h$.

\subsection{Solitary wave running over a flat bottom}
\begin{figure}[!t]
\center
\includegraphics[width=4.75in, trim={0 -0.25in 0 -0.25in}, clip, frame]{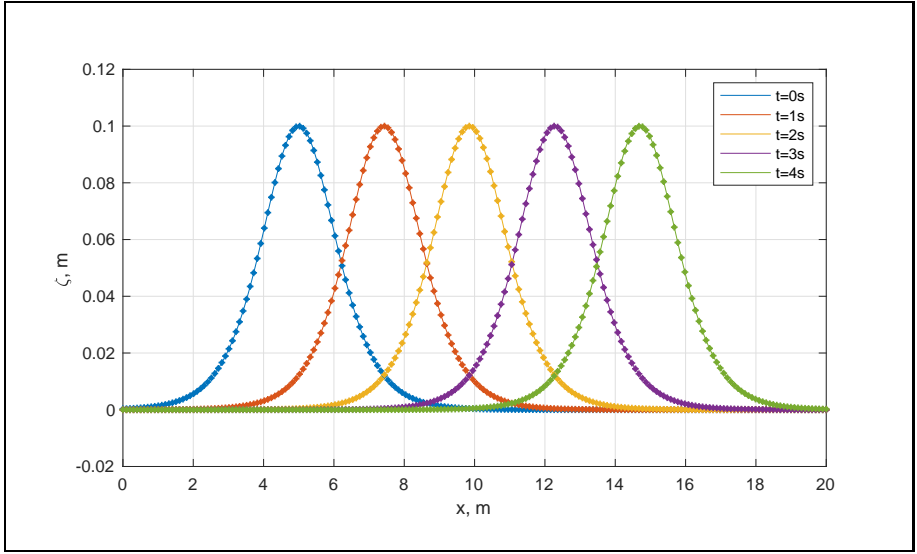}
\caption{Free surface elevations at specified time points for a solitary wave propagation over a flat bottom obtained numerically (solid line) and analytically (diamonds).}
\label{Fig:E1}
\end{figure}

In this example a numerical simulation of a solitary wave propagating over a flat bottom is compared with the analytical solution for this process that satisfies the Green-Naghdi equations in Eq.(\ref{Eq:GN}) for $\alpha=1$ \cite{duran_and_marche_2017}:
\begin{linenomath}
\begin{equation}
h(x,t) = H_0 + a_0\sech^2\left(\kappa(x-x_0-c_0t)\right), 
\quad
h\mathbf u(x,t) = c_0 h(x,t) - c_0 H_0,
\label{Eq:SW}
\end{equation}
\end{linenomath}
where $H_0$ is the reference water level, $a_0$ is the solitary wave height, $x_0$ is the initial wave position, and
\begin{linenomath}
\begin{equation}
\kappa = \frac {\sqrt{3 a_0}} {2H_0 \sqrt{H_0+a_0}},
\quad
c_0 = \sqrt{g (H_0+a_0)}.
\end{equation}
\end{linenomath}

A problem domain $\Omega = (0,20)\times(-2.5\cdot10^{-2},2.5\cdot10^{-2})\,m^2$ is partitioned into a finite element mesh comprised of $400\times1$ square cells containing 2 triangular elements. The parameters for the solitary wave in this experiment are set as: $H_0=0.5m$, $a_0=0.1m$, and $x_0=5m$. All of the boundaries of the domain are set as reflecting boundaries while the initial position of the solitary wave ensures there is negligible interaction between the solitary wave and the boundaries at $x=0$ and $x=20m$.  The initial conditions are set with Eq.(\ref{Eq:SW}) at $t=0$, and a second-order Runge-Kutta time stepping is used to propagate the numerical solution with a time step $\Delta t=10^{-4}s$. Results of the numerical simulation are presented in Fig.\ref{Fig:E1} where the free surface elevation solution is compared with the analytical solution at $t=\{0,1,2,3,4\}s$. The results suggest that the numerical solution closely matches the analytical solution for the solitary wave propagation over a flat bottom.

\subsection{Shoaling solitary wave reflected off a solid wall}
\begin{figure}[!t]
\center
\includegraphics[width=4.75in, trim={0 0 0 0}, clip, frame]{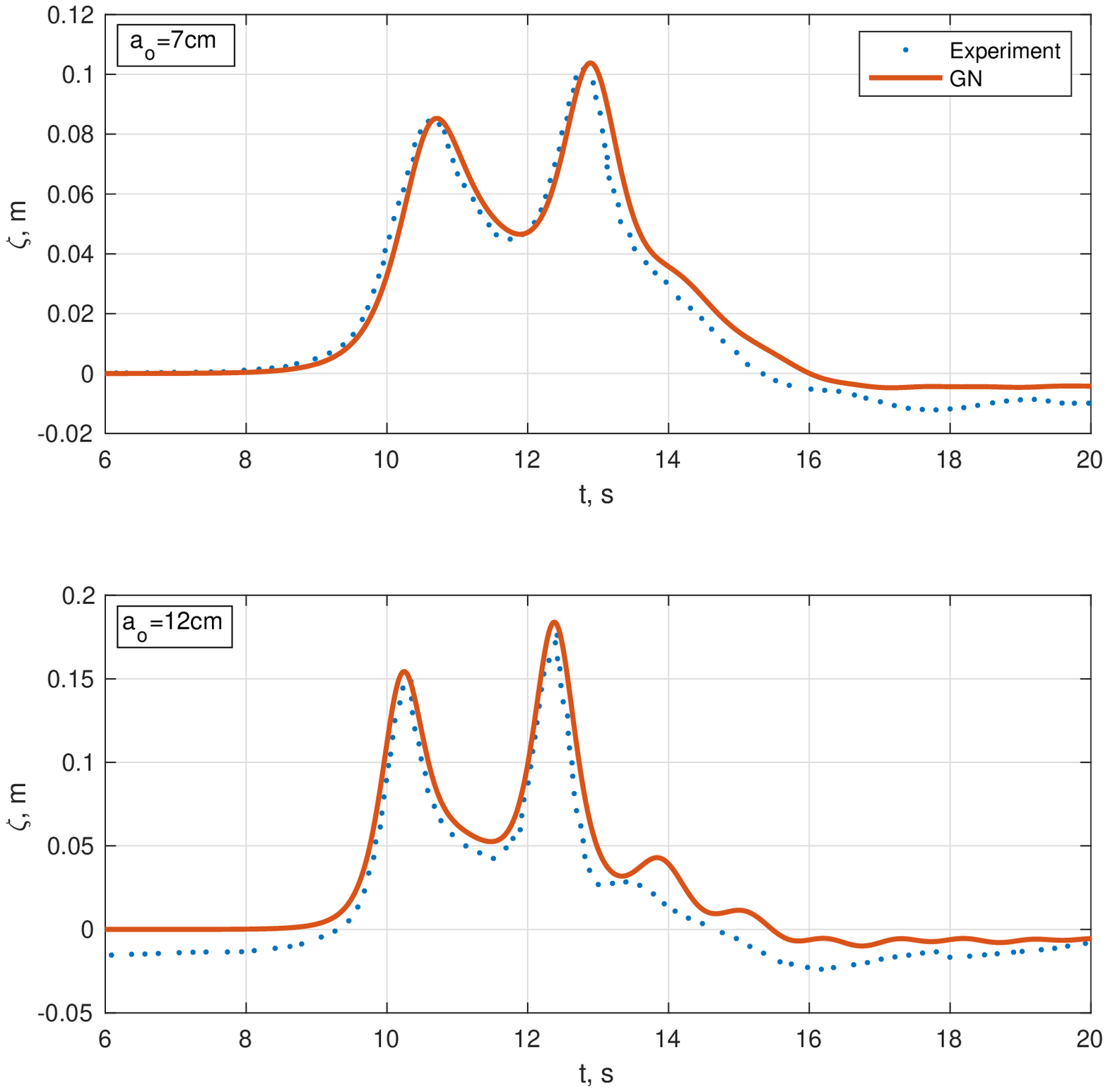}
\caption{Free surface elevations at $x=17.75m$ for the Green-Naghdi (GN) equations simulations and experimental results reported by Walkley and Berzins in \cite{walkley_and_berzins_1999}.}
\label{Fig:E2}
\end{figure}

This numerical experiment simulates solitary waves shoaling over a linearly sloping bottom and reflecting off a solid wall. Results of the numerical experiment are compared with the physical experiment results obtained by Dodd in \cite{dodd_1998} as reported by Walkley and Berzins in \cite{walkley_and_berzins_1999}. Two solitary waves from Eq.(\ref{Eq:SW}) have been studied in the physical experiment with the following parameters: $H_0=0.7m$, $a_0=\{0.07, 0.12\}m$.

For this numerical experiment, a finite element mesh of a problem domain $\Omega = (-20,20)\times(-1.25\cdot10^{-1},1.25\cdot10^{-1})\,m^2$ is generated with $160\times1$ square cells containing 2 triangular elements. The domain has a flat bottom that starts sloping up at $1:50$ rate starting from $x=0$. The solitary waves in this experiment are initially positioned at $x_0=-10m$, the initial conditions are set with Eq.(\ref{Eq:SW}) at $t=0$, and a second-order Runge-Kutta time stepping is used to propagate the numerical solution for $20s$ with a time step $\Delta t=10^{-3}s$. Similarly to the previous experiment, all of the boundaries of the domain are set as reflecting boundaries while setting the initial positions of the solitary waves to ensure negligible interaction with the boundary at $x=-20m$.  The results of the numerical simulations are presented in Fig.\ref{Fig:E2} where the free surface elevation measurements at $x=17.75m$ are compared with the experiment results reported by Walkley and Berzins. The numerical simulation results are in good agreement with the experimental results, the hydrodynamic model is able to capture the peak free surface elevations as the waves pass the measurement point before and after the reflection. 

\subsection{Regular waves running over a trapezoidal bar}
\begin{figure}[!t]
\center
\includegraphics[width=4.75in, trim={0 0 0 0}, clip, frame]{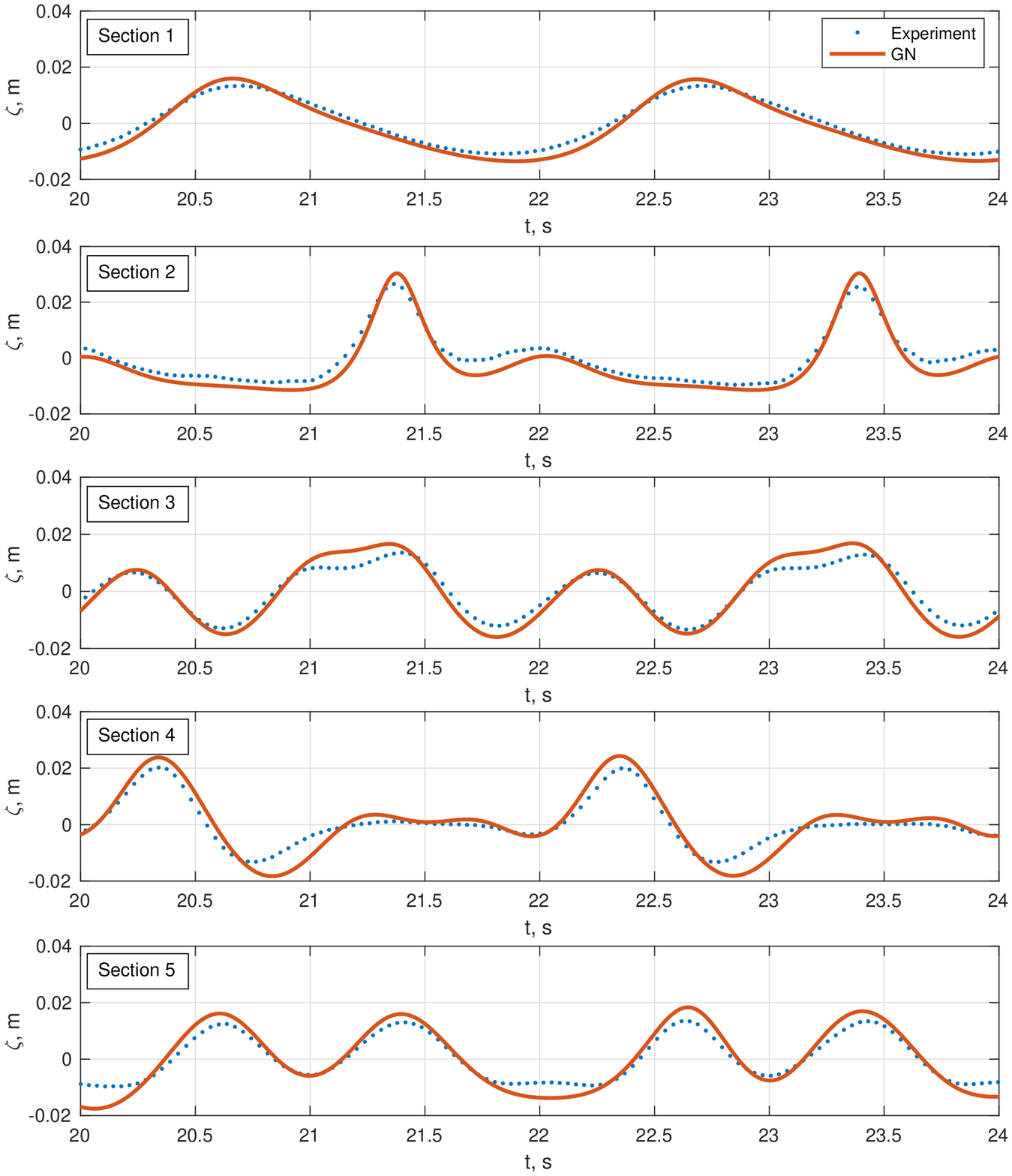}
\caption{Free surface elevation measurement at 5 sections obtained with the Green-Naghdi equations simulation and compared with the experimental results by Dingemans in \cite{dingemans_1994}.}
\label{Fig:E3}
\end{figure}

This numerical experiment simulates the propagation of regular waves over a fully submerged trapezoidal bar. The bathymetry for this experiment is defined as
\begin{linenomath}
\begin{equation}
b(x)=
\begin{cases}
-0.4+\frac{x-6}{20}&\text{if}\,\,\,\,6 < x < 12 \\
-0.1&\text{if}\,\,\,\,12\leq x < 14 \\
-0.1-\frac{x-14}{10}&\text{if}\,\,\,\,14\leq x < 17 \\
-0.4&\text{elsewhere}
\end{cases}.
\end{equation}
\end{linenomath}
The water is initially at rest, and a regular wave with the amplitude of $0.01m$ and the period of $2.2s$ is generated at the left boundary of the domain. The nontrivial bathymetry in this experiment induces substantial non-linearities: while shoaling and steepening over the upward slope, the regular waves generate higher-harmonics that are progressively released over the downward slope \cite{duran_and_marche_2017}.

In this numerical simulation, a finite element mesh of a problem domain $\Omega = (0,30)\times(-2.5\cdot10^{-2},2.5\cdot10^{-2})\,m^2$ is generated with $600\times1$ square cells containing 2 triangular elements. The left boundary of the domain functions as a wavemaker boundary where free-surface elevation is prescribed to generate regular waves. The rest of the boundaries are set as reflecting boundaries. A second-order Runge-Kutta time stepping is used to propagate the numerical solution with the time step $\Delta t=5\cdot10^{-3}s$. The free surface elevation measurements were performed at 5 different sections located at $x=\{10.5,13.5,15.7,17.3,19\}m$. The results of these measurements along with their comparison to the physical experiment results by Dingemans in \cite{dingemans_1994} are presented in Fig.\ref{Fig:E3}. The results suggest that the dispersive wave hydrodynamic model is able to capture the regular wave transformations over the trapezoidal bar with a sufficient accuracy both in the wave amplitude and the frequency.

\subsection{A solitary wave over a sloping beach with mobile bed}
\begin{figure}
\center
\includegraphics[width=4.75in, trim={0 0.7in 0 0}, clip, frame]{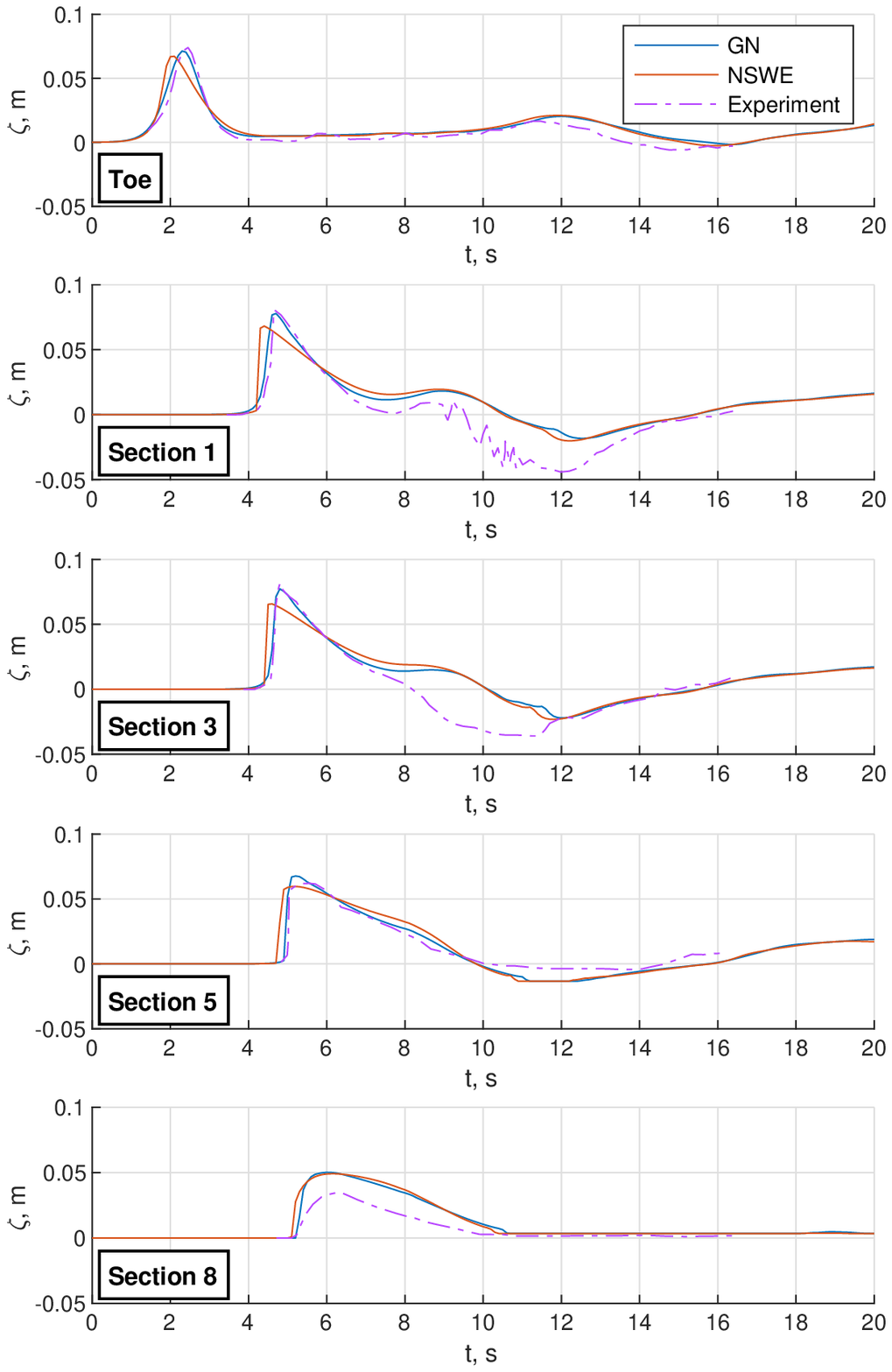}
\caption{The free surface elevation measurements at stations for the Green-Naghdi (GN) and nonlinear shallow water equations (NSWE) simulations, and experimental results by Sumer \emph{et al.} \cite{sumer_etal_2011}.}
\label{Fig:Wave}
\end{figure}

The model has been validated against the experiment conducted by Sumer \emph{et al.} \cite{sumer_etal_2011} to measure flow and bed morphology induced by a solitary wave over a sloping beach. In the experiment four solitary waves have been run over a sloping beach inclined at $1:14$ rate; and, subsequently, a number of measurements have been performed, such as extents of sediment erosion and deposition over the sloping beach, and the free surface elevation at nine measuring stations. The measuring stations are located at the toe of the sloping beach, and at eight sections located 4.63, 4.69, 4.87, 5.11, 5.35, 5.59, 5.65, and 5.85 meters from the toe. 

The choice of this experiment for the validation of the model has been motivated by the following reasons: (1) dispersive wave effects are prevalent in this experiment, and the Green-Naghdi equations should be used to resolve accurately the water wave dynamics, (2) in this experiment the solitary waves have sufficiently high amplitude to experience wave breaking; therefore, a wave breaking detection is required to switch to the nonlinear shallow water equations in surf zones, (3)  in this experiment the sloping beach undergoes substantial sediment erosion and deposition that affect the bed surface elevation of the beach. Thus, performing numerical simulations of this experiment and comparing the results to the experimental ones have the potential to showcase all key features of the presented numerical model, such as the ability of the Green-Naghdi equations to simulate accurately water motion and capture dispersive wave effects, capacity of the numerical model to detect wave breaking regions and switch to the nonlinear shallow water equations in such regions, and the facility of the model to estimate sediment erosion and deposition due to bed-load transport.     

To carry out the numerical simulations, a problem domain $\Omega = (-10,10)\times(-2.5\cdot10^{-2},2.5\cdot10^{-2})\,m^2$ is partitioned into a finite element mesh comprised of $400\times1$ square cells containing 2 triangular elements. The sloping beach toe is located at $x=0$, and all boundaries of the mesh are specified as reflecting boundaries which reflects the physical experiment setup. A two stage second-order Runge-Kutta method is used to perform time integration with a time step $\Delta t=5\cdot10^{-3}s$. The initial conditions for this experiment are characterized by solitary waves from Eq.(\ref{Eq:SW}) at $t=0$ with the reference water level $H_0=0.4m$, the solitary wave height $a_0=0.071m$, the initial wave position $x_0=-5m$. Finally, for these simulations the bottom friction force is introduced into the numerical model through the source term $\boldsymbol S(\boldsymbol q)$ by setting
\begin{linenomath}
\begin{equation}
\mathbf f = C_f{\vert\mathbf{u}\vert\mathbf{u}},
\label{Eq:f}
\end{equation}
\end{linenomath}
where the Chezy friction coefficient $C_f=0.012$.

\begin{figure}[t]
\center
\includegraphics[width=4.75in, frame]{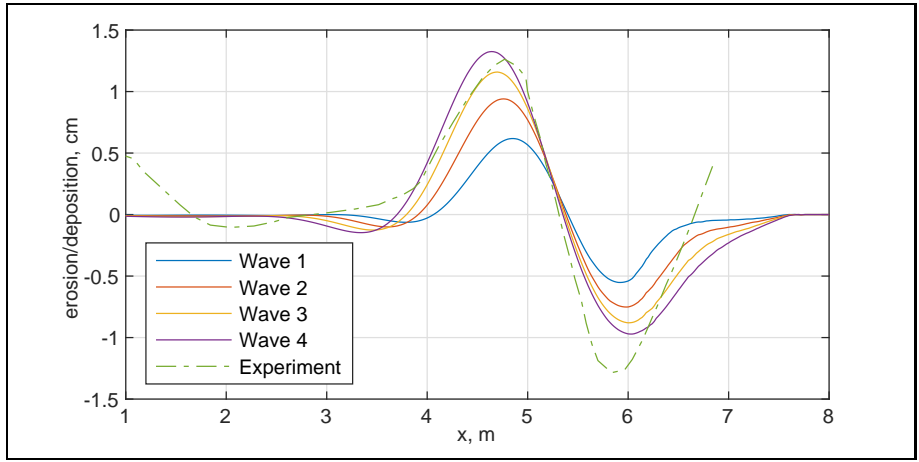}
\caption{Sediment erosion and deposition over the sloping beach for the decoupled model approach, and the experiment by Sumer \emph{et al.} \cite{sumer_etal_2011}.}
\label{Fig:Erosion1}
\end{figure}

\begin{figure}[t]
\center
\includegraphics[width=4.75in, frame]{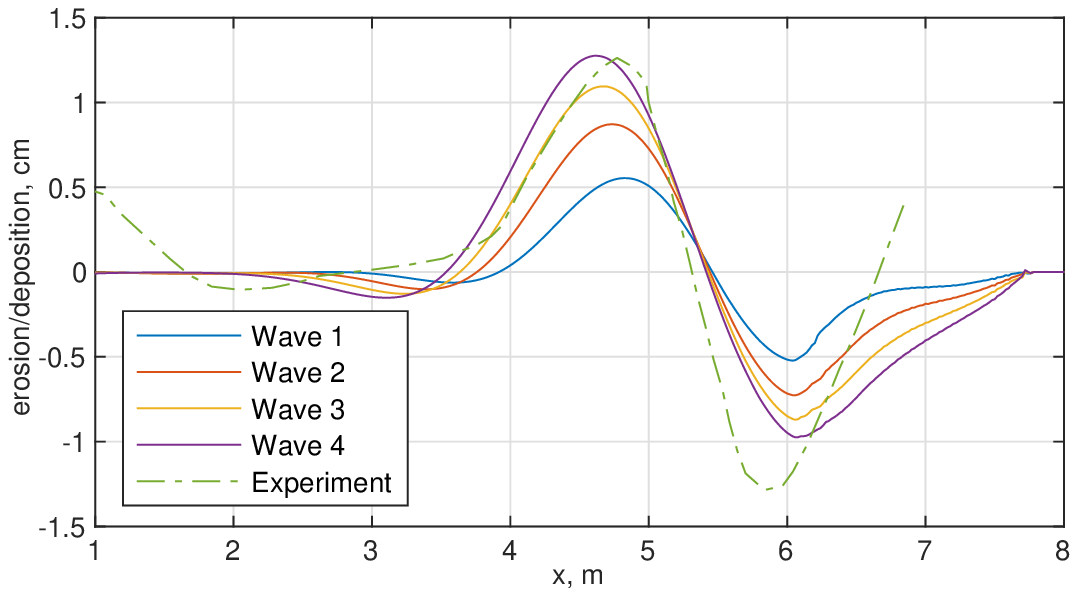}
\caption{Sediment erosion and deposition over the sloping beach for the coupled model approach, and the experiment by Sumer \emph{et al.} \cite{sumer_etal_2011}.}
\label{Fig:Erosion2}
\end{figure}

The simulations have first been  performed over a rigid bed to validate the hydrodynamic model. Using separately the Green-Naghdi and nonlinear shallow water equations, the simulations have been run for 20 seconds which is a sufficient time for solitary waves to run up and run down along the sloping beach in this experiment. Fig.\ref{Fig:Wave} presents the free surface elevations obtained at the measuring stations from the experiment by Sumer \emph{et al.} \cite{sumer_etal_2011}. As expected, in terms of accuracy the Green-Naghdi equations substantially outperform the nonlinear shallow water equations in the run up stage at the measuring stations located offshore. It is also evident that solitary waves break too early in the nonlinear shallow water equations simulations. In fact, the experimental results suggest that wave breaking occurs somewhere between the sections 3 and 5 which is accurately captured by the Green-Naghdi equations. However, neither model is able to accurately capture the water motion in the swash zone as evidenced by the free surface elevation measurements at the onshore section 8. We believe that these inaccuracies are due to the nontrivial physics that govern the water motion in swash zones, and to the limitations of the wetting-drying algorithm used in the simulations. Subsequently, the models are unable to capture correctly the water motion during the run down stage of the simulations. Nonetheless, the results are deemed satisfactory given the complexity of the physical processes occurring in flows induced by solitary waves over a sloping beach.      

The full dispersive wave hydro-morphodynamic model has been used for the erodible bed simulations. The Grass model \cite{grass_1981} in Eq.(\ref{Eq:Qb}) for the sediment flux $\mathbf{Q}_b$ has been used with $A=4.75\cdot10^{-3}$. Both the decoupled and coupled models have been used in the simulations. In each simulation a total of 4 solitary waves have been run over the sloping beach such that the initial conditions except the bathymetry are reset before each wave run. Each run has been performed for 2 minutes and 30 seconds which is a sufficient time for water to substantially settle after the run down. In this experiment the time scales in the hydrodynamic and morphodynamic models are comparable \cite{li_etal_2019}. Therefore, in the decoupled model the bed surface update has to be performed every time step of the hydrodynamic model. The bed surface erosion and deposition results obtained from the simulations are presented in Fig.\ref{Fig:Erosion1} for the decoupled model and in Fig.\ref{Fig:Erosion2} for the coupled model. As expected, the bed surface evolution in the offshore area is accurately estimated by both models since the hydrodynamic part captures the water motion in that area with sufficient accuracy. On the other hand, in the onshore area the models capture sediment erosion and deposition less accurately due, in part, to low accuracy of the hydrodynamic model in the swash zone. Overall, the results are considered satisfactory and indicate a promise for further development of the presented hydro-morphodynamic model, e.g. towards the extension of the model with suspended-load transport. Moreover, the decoupled model performed well relative to the coupled model and can provide a viable alternative to the coupled model, in particular, in cases where the time scales in the hydrodynamic part are shorter than in the morphodynamic part.         

\subsection{The Faro-Olh\~ao inlet}
\begin{figure}
\centering
\includegraphics[width=4.75in]{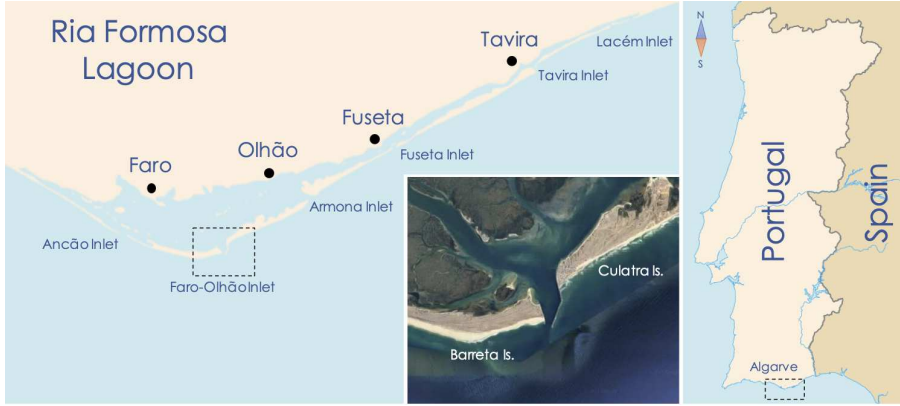}
\caption{The Ria Formosa lagoon diagram.}
\label{Fig:Ria}
\end{figure}

The presented work is a stepping stone in an ongoing project aimed at modeling hydro-morphodynamic processes in the Ria Formosa lagoon in the Algarve Region of Portugal (cf. Fig.\ref{Fig:Ria}). Although modeling hydro-morphodynamic processes in the lagoon are beyond the scope of this paper, the hydrodynamic model has been used to simulate water waves around the Faro-Olh\~ao inlet of the lagoon. The lagoon stretches about 55km along the southern coast of Portugal. It is separated from the Atlantic Ocean by a series of barrier islands, and has six naturally occurring and artificial tidal inlets. Astronomical tides in the area lead to nearly 2m in water level variation \cite{pacheco_etal_2007}, which has formed large salt marshes and mudflats in the lagoon. The salt marshes and mudflats cover nearly 70\% of the total area of the lagoon \cite{carrasco_etal_2018}. The lagoon is a valuable regional resource for tourism and fisheries, and a natural habitat for various protected species. The Faro-Olh\~ao inlet provides the gateway for the channels that connect two main cities of the region, Faro and Olh\~ao, to the open ocean. The inlet was artificially opened and stabilized with jetties between 1929 and 1955 \cite{pacheco_etal_2007}. 

Four sources of the region bathymetric data have been identified: Portuguese Hydrographic Institute bathymetric model \cite{HIP}, bathymetric surveys performed under SCORE project \cite{SCORE, gorbena_etal_2017}, LiDAR bathymetric data of the coast of Portugal \cite{gorbena_etal_2017, SCORE_DB}, and EMODnet Bathymetry Digital Terrain Model (DTM 2018) data \cite{EMODnet}. The data sources have varying levels of detail and coverage: (1) EMODnet is the single source that has bathymetric data for the open ocean, but it is less detailed in the near shore areas and inside the lagoon; (2) Portuguese Hydrographic Institute data has great resolution in the lagoon channels only; (3) LiDAR data has the best resolution within the lagoon, but it is missing bathymetric data where the water is too deep to perform LiDAR measurements; (4) SCORE project data has a good representation of the lagoon, but it is missing some features, such as jetties at Faro-Olh\~ao and Tavira inlets. The data has been combined to generate a single bathymetric profile of the lagoon and its surrounding area.

The resulting bathymetric profile has been used to generate a finite element mesh, which has over $10^5$ triangular elements with diameters ranging between 10m and 200m, of the Faro-Olh\~ao inlet and its surrounding area. Intuitively, most of the flow variation occurs in and around the inlet where the mesh has a finer resolution. A second-order Runge-Kutta time stepping is used for time integration with a time step $\Delta t=10^{-1}s$. The waves in the model are forced at the open ocean boundary of the problem domain with M2 tidal constituent with the amplitude of 1.01m and the period of 12.42h \cite{dias_etal_2009}. The dispersive correction is not applied in the vicinity of the open ocean boundary, and the tidal wave is imposed through boundary conditions for the nonlinear shallow water equations. The bottom friction force is introduced into these numerical simulations through setting $\mathbf f$ in the source term as in Eq.(\ref{Eq:f}) with the Chezy friction coefficient $C_f=0.0045$.

\begin{figure}
\centering
\begin{subfigure}{0.475\linewidth}
\includegraphics[trim={2.5cm 4.5cm 2cm 4.5cm},clip,angle=270,origin=c,width=2.2in,frame]{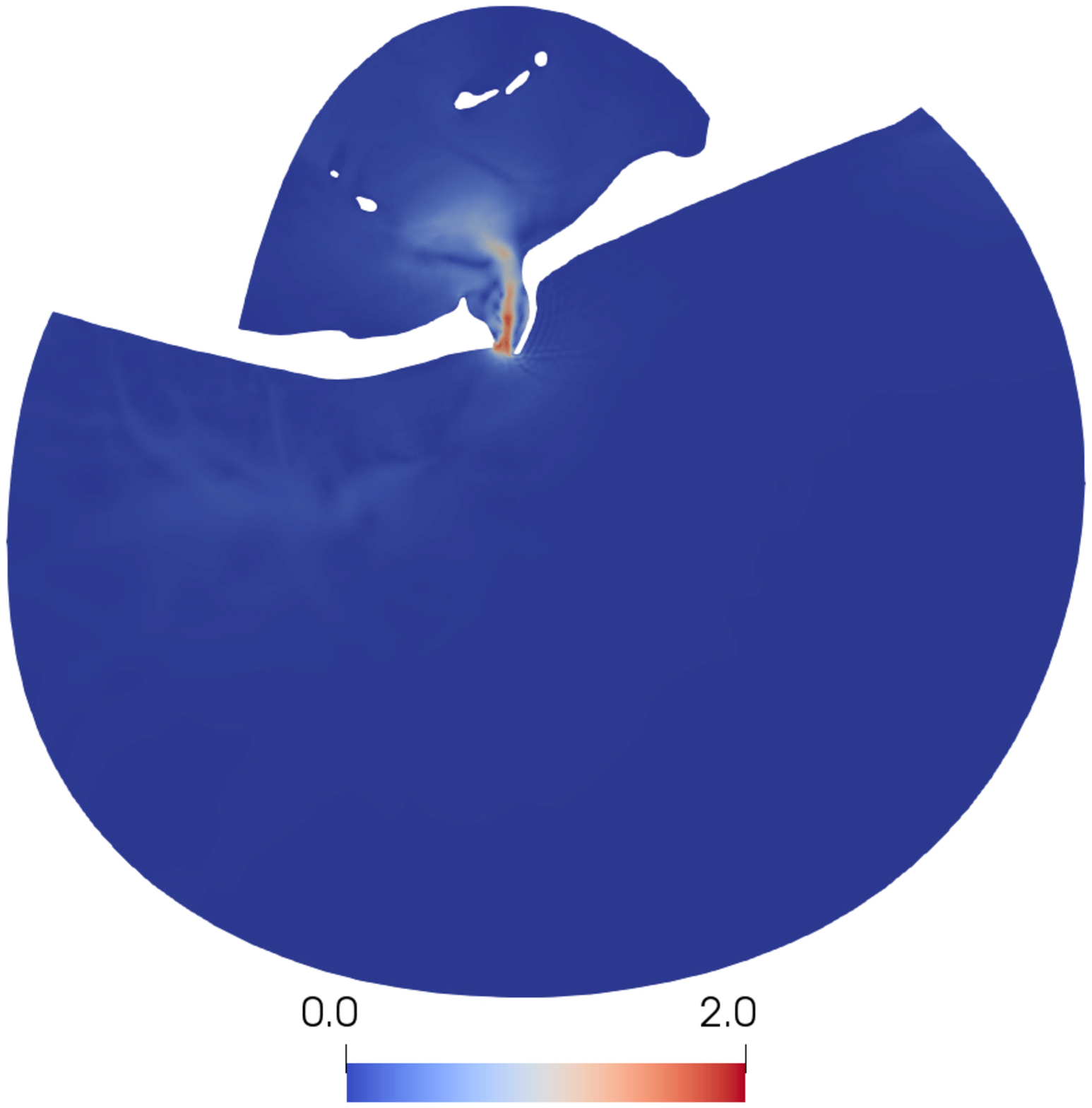}
\caption{$\lvert\textbf{u}\rvert_{\text{GN}}$}
\end{subfigure}
\begin{subfigure}{0.475\linewidth}
\includegraphics[trim={2.5cm 4.5cm 2cm 4.5cm},clip,angle=270,origin=c,width=2.2in,frame]{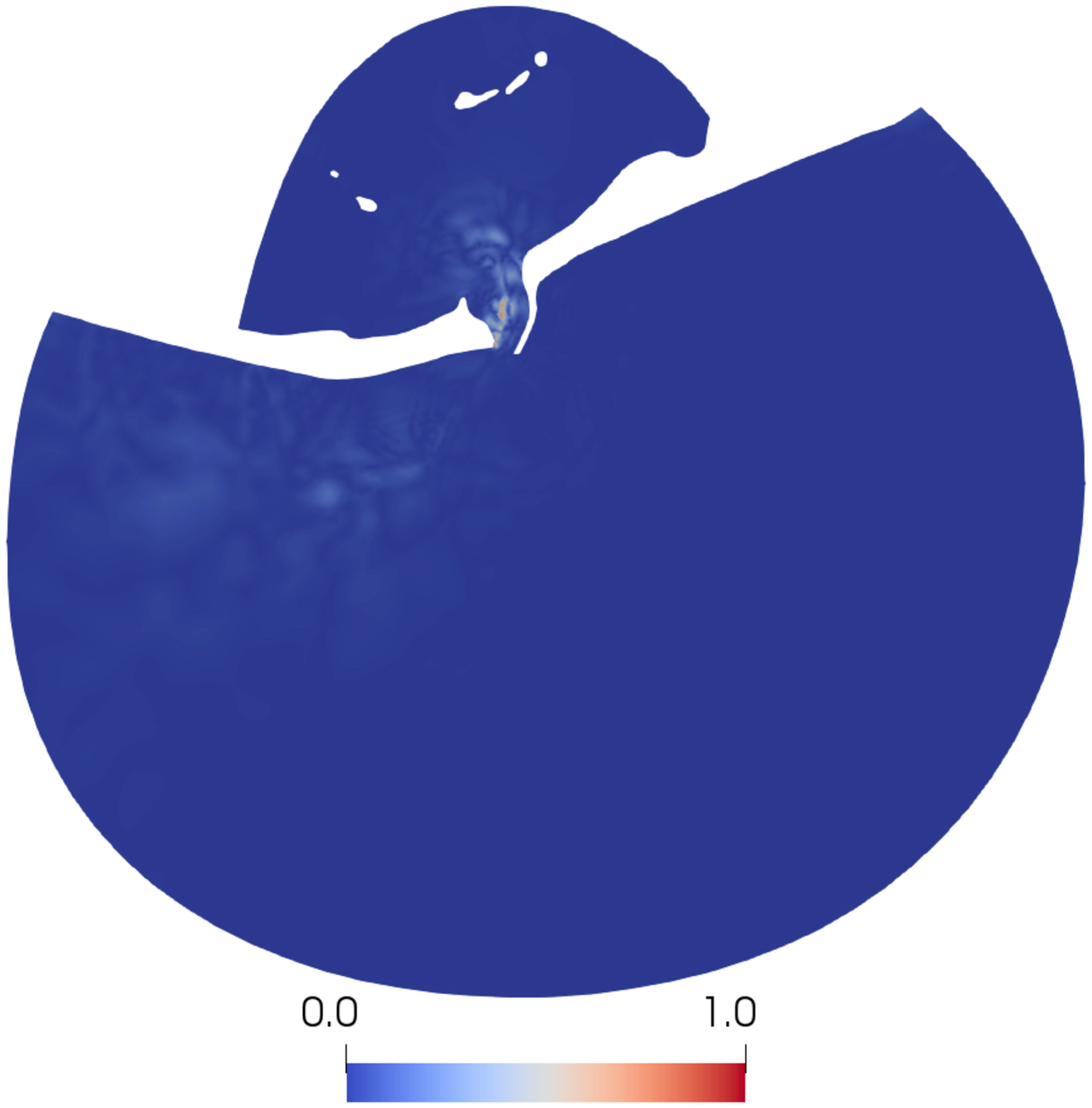}
\caption{$\lvert\lvert\textbf{u}\rvert_{\text{GN}}-\lvert\textbf{u}\rvert_{\text{NSWE}}\rvert$}
\end{subfigure}
\caption{Velocity fields around the time of the peak inflow velocity at the Faro-Olh\~ao inlet.}
\label{Fig:maxIN}
\end{figure}

\begin{figure}
\centering
\begin{subfigure}{0.475\linewidth}
\includegraphics[trim={2.5cm 4.5cm 2cm 4.5cm},clip,angle=270,origin=c,width=2.2in,frame]{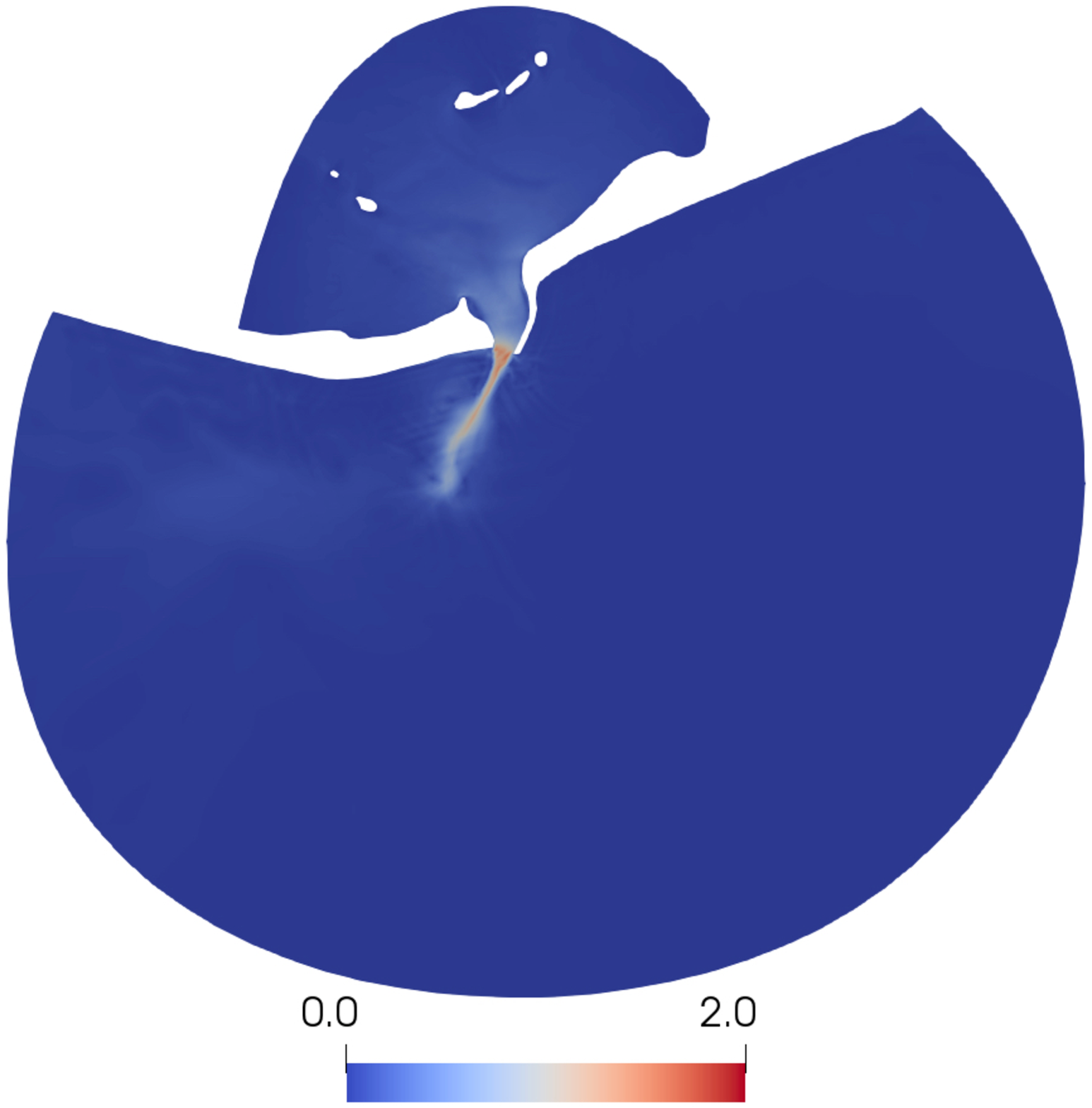}
\caption{$\lvert\textbf{u}\rvert_{\text{GN}}$}
\end{subfigure}
\begin{subfigure}{0.475\linewidth}
\includegraphics[trim={2.5cm 4.5cm 2cm 4.5cm},clip,angle=270,origin=c,width=2.2in,frame]{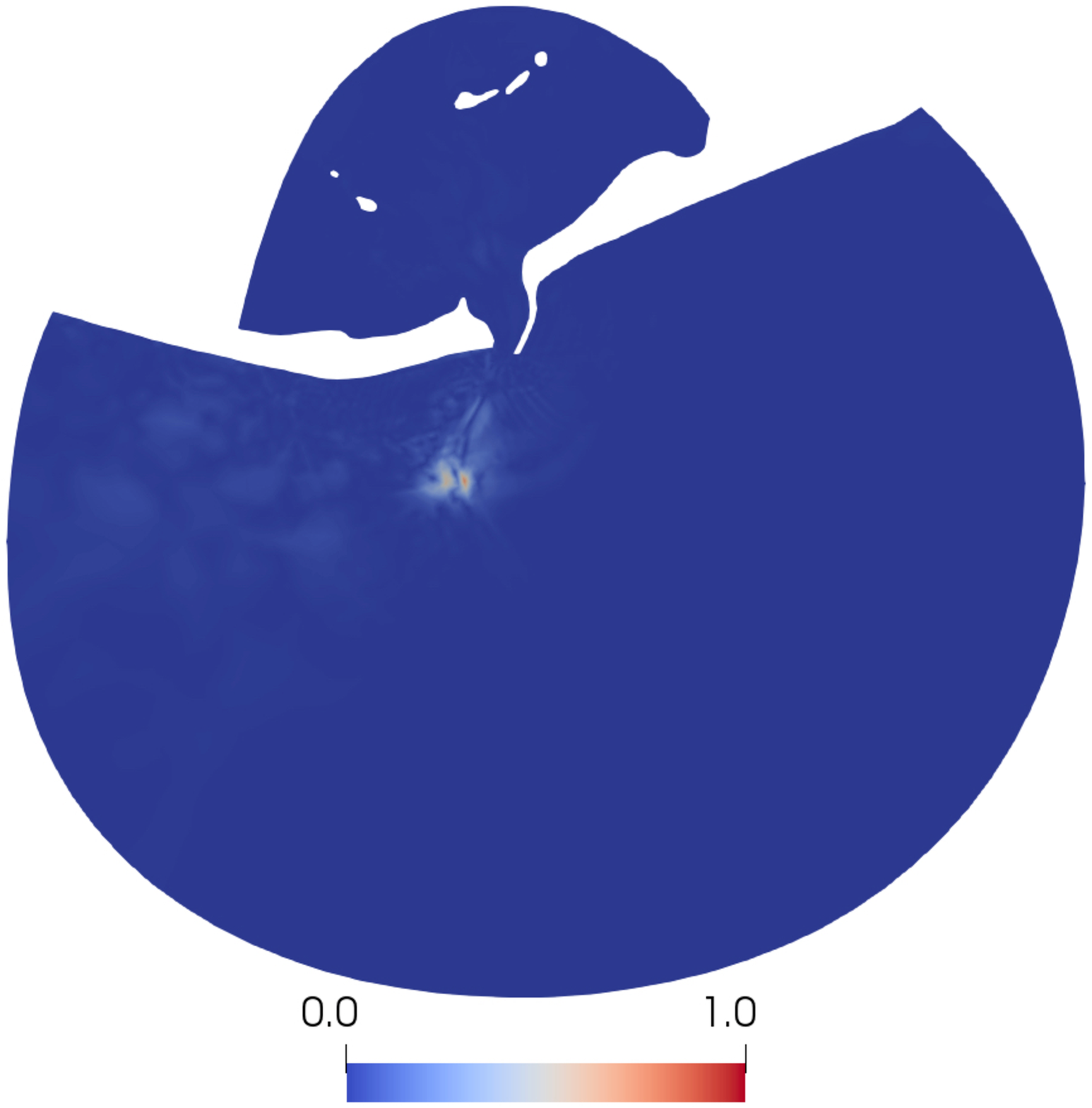}
\caption{$\lvert\lvert\textbf{u}\rvert_{\text{GN}}-\lvert\textbf{u}\rvert_{\text{NSWE}}\rvert$}
\end{subfigure}
\caption{Velocity fields around the time of the peak outflow velocity at the Faro-Olh\~ao inlet.}
\label{Fig:maxOUT}
\end{figure}

The model has been used to simulate water waves for 2 days with the Green-Naghdi and nonlinear shallow water equations. The velocity profiles around the times of the peak inflow and outflow velocities at the neck of the inlet are presented in Fig.\ref{Fig:maxIN} and Fig.\ref{Fig:maxOUT}, respectively. The hydrodynamic model is able to successfully simulate the water waves with the Green-Naghdi equations over the irregular shaped unstructured mesh. The magnitude of difference between the velocity profiles obtained with the Green-Naghdi and nonlinear shallow water equations shows that there is a considerable dissimilarity between these two computations. The results indicate the potential of using the model to simulate water waves in coastal regions with irregular geometries, and,  subsequently, to study morphodynamic processes.       

\section{Conclusions}
In this paper a hydro-morphodynamic model that couples a depth-averaged dispersive water wave model, the Green-Naghdi equations, with the Exner equation has been introduced. Although there are numerous works that couple the nonlinear shallow water equations with the Exner equation, to the best of authors' knowledge, the coupling of the sediment continuity model with the Green-Naghdi equations has not been attempted before this work. The presented model is well suited for studying the bed surface evolution under bed-load transport in areas where dispersive wave effects are prevalent and should thus be included in the hydrodynamic model. 

Numerical methods that utilise discontinuous Galerkin finite element methods have been presented for the hydro-morphodynamic model. The Strang operator splitting technique has been employed to single out the dispersive part of the Green-Naghdi equations for separate treatment. The resulting numerical models are augmented with wetting-drying, breaking wave detection, and slope limiting features. The numerical solution algorithm developed for the hydrodynamic part of the coupled model has been validated with three numerical experiments that have displayed the capacity of the algorithm to accurately capture hydrodynamics of both solitary and regular waves. The full hydro-morphodynamic model has also been used to simulate flow and sediment transport induced by solitary waves over a sloping beach. Comparing the numerical results with the experimental results collected by Sumer \emph{et al.} \cite{sumer_etal_2011}, the numerical experiments have demonstrated that the presented model is capable of modeling water waves and sediment transport with a satisfactory accuracy in areas where the wave dispersion effects prevail. Moreover, the dispersive wave hydrodynamic model has been used to simulate water waves in the Faro-Olh\~ao inlet of the Ria Formosa lagoon in Portugal. The results indicate that the model is capable of performing simulations over irregular shaped unstructured meshes. This capability is important for simulations required in studies of hydro-morphodynamic processes in coastal areas.  

The presented hydro-morphodynamic model has a potential to be used in simulations of hydro-morphodynamic processes caused by dispersive waves in large coastal areas. While the hydrodynamic part of the model is capable of capturing water wave dynamics with a sufficient accuracy up to swash zones, the properly calibrated morphodynamic part of the model can estimate bed erosion/deposition due to bed-load transport. As a further development, the presented model can be extended with suspended-load transport to enhance the ability of the presented model to estimate morphological changes of coastal areas induced by the actions of waves and currents.       

\section{Acknowledgments}
This work has been supported by funding from the National Science Foundation Grant 1854986, and the Portuguese government through Funda\c{c}\~ao para a Ci\^encia e a Tecnologia (FCT), I.P., under the project DGCOAST (UTAP-EXPL/MAT/0017/2017). The authors of the paper would like to acknowledge the support of the Texas Advanced Computing Center through the allocation TG-DMS080016N used in the parallel computations of this work. The authors acknowledge the support of researchers at the Centre for Marine and Environmental Research of the University of Algarve in obtaining necessary data for the Ria Formosa lagoon and its surrounding region. The authors thank anonymous reviewers of this paper for their comments and suggestions that helped improve the paper.

\bibliography{mybibfile}

\end{document}